\documentclass[12pt]{article}
\usepackage{eurosym}
\usepackage[left=2.2cm, right=2.2cm, top=2.2cm, bottom=2.2cm]{geometry}
\usepackage{amsfonts}
\usepackage{amssymb}
\usepackage{amsmath}
\usepackage{amsthm}
\usepackage{lscape}
\usepackage{graphicx}
\usepackage{pgfplots}
\usepackage{tikz}
\usepackage{epsfig}
\pgfplotsset{compat=newest}
\pgfplotsset{plot coordinates/math parser=false}
\pgfplotsset{ignore zero/.style={%
		#1ticklabel={\ifdim\tick pt=0pt \else\pgfmathprintnumber{\tick}\fi}
}}
\newlength\figureheight
\newlength\figurewidth

\DeclareMathOperator*{\plim}{plim}

\newcommand{\Keywords}[1]{\par\noindent{{\em \large{Keywords}\/}: #1}}

\def\1{1\!{\rm l}}
\newcommand{\Plim}[1]{\raisebox{0.5ex}{\scalebox{0.8}{$\displaystyle \plim_{#1}\;$}}}

\begin{document}

\title{Indirect Inference With(Out) Constraints}
\author{David T. Frazier\thanks{Department of Econometrics and Business Statistics, Monash University,
Melbourne, Australia.} 
and Eric Renault\thanks{Department of Economics, University of Warwick.}\footnote{We thank the Co-Editor, Andres Santos, and three anonymous referees for many helpful comments that have greatly improved the paper. In addition, we thank Geert Dhaene, Philipp Ketz, and Mervyn Silvapulle for helpful discussions. }}
\maketitle

\begin{abstract}
\noindent Indirect Inference (I-I) estimation of structural parameters $\theta$ {{requires matching observed and simulated statistics, which are most often generated using an auxiliary model that depends on instrumental parameters $\beta$.}} {The estimators of the instrumental parameters will encapsulate} the statistical information used for inference about the structural parameters. As such, artificially constraining these parameters may restrict the ability of the auxiliary model to accurately replicate features in the structural data, which may lead to a range of issues, such as, a loss of identification. However, in certain situations the parameters $\beta$ naturally come with a set of $q$ restrictions. Examples include settings where $\beta$ must be estimated subject to $q$ possibly strict inequality constraints $g(\beta) > 0$, such as, when I-I is based on GARCH auxiliary models. In these settings we propose a novel I-I approach that uses appropriately modified unconstrained auxiliary statistics, which are simple to compute and always exists. We state the relevant asymptotic theory for this I-I approach without constraints and show that it can be reinterpreted as a standard implementation of I-I through a properly modified binding function. Several examples that have featured in the literature illustrate our approach. 
\end{abstract}
\hspace{6cm}


\section{Introduction}

The indirect estimation procedures of Gourieroux, Monfort and Renault (1993)
(hereafter, GMR), Smith (1993) and Gallant and Tauchen (1996) (hereafter, GT) provide convenient
estimation methods when efficient estimation of a fully parametric
structural model is a daunting task due to the intractability of the
likelihood function. GMR motivate Indirect Inference (I-I) by arguing that in such cases a natural procedure is to replace the
likelihood function by another criterion based on some convenient auxiliary
(or naive) model that is simpler but possibly misspecified. The overall aim of I-I is then to conduct correct inference ``based on this incorrect criterion."

As described by Jiang and Turnbull (2004), the ``essential
ingredients" of I-I are as follows:

\noindent(i) A parametric model for data generation, with distribution $P_{\theta }$
that depends on an unknown vector $\theta \in \Theta \subset 
\mathbb{R}
^{d_{\theta }}$ of parameters of interest. This model is the so-called
structural model and $\theta $ is the vector of structural parameters.

\noindent(ii) An intermediate or auxiliary statistic, say $\hat{\beta}_{T}$, of dimension $d_{\beta }\geq d_{\theta }$, which is a functional of the
observed sample $\left\{ y_{t}\right\} _{t=1}^{T}.$

\noindent(iii) A bridge (or binding) relationship $\beta =b(\theta )$ defined
between the true unknown value $\theta ^{0}$ of the structural parameters
and $\beta ^{0}=\Plim{T\rightarrow\infty} \hat{\beta}_{T}$, where the unknown quantity $\beta
^{0}=b(\theta ^{0})$ is referred to as the pseudo-true value of the auxiliary
parameters.

\noindent(iv) With the auxiliary statistic $\hat{\beta}_{T}$ replacing $\beta $, the
bridge relationship is used to compute an I-I estimator of $\theta$ by ``inverting'' $b(\theta)$. 

Jiang and \ Turnbull (2004) acknowledge that ``the choice of an intermediate
statistic $\hat{\beta}_{T}$ is not necessarily unique'', however, the authors argue that ``in any given
situation there is often a natural one to use." Herein, we question this traditional interpretation of I-I when it pertains to examples where the
choice of ``intermediate statistic'' is
ambiguous {due to the fact that the parameters of the {auxiliary model used in I-I}} must be estimated subject to a vector of inequality constraints.  

In this commonly encountered situation, the choice of appropriate intermediate statistics for I-I can be ambiguous for several reasons: firstly, as noted
by Calzolari, Fiorentini and Sentana (2004) (hereafter, CFS), the pseudo-likelihood
function of the auxiliary model may not be well-behaved when certain
parameter restrictions are violated and, hence, {without these additional restrictions $\hat{\beta}_{T}$} cannot be obtained; secondly, if the pseudo-true value $\beta^0$ is on (or near) the boundary of the parameter space defined by the inequality constraints, the intermediate statistic $\hat{\beta}_T$ may be insufficient to identify $\theta^0$; lastly, even if identification of $\theta^0$ is possible, if $\beta^0$ is on (or near) the boundary of the parameter space defined by the constraints, pseudo-maximum likelihood (hereafter, PML) estimation will lead to an
intermediate statistic that is not well-suited for I-I because it is not asymptotically normal (see Andrews, 1999 and CFS for details). 

The question then is how to choose the intermediate statistics so as to guarantee consistent and asymptotically normal I-I estimators of $\theta^0$ even though $\beta^0$  can lie on (or near) the boundary of the parameter space defined by these inequality constraints? One approach, which is proposed by CFS, is to consider {as our intermediate statistic}  a ``well-behaved'' linear combination of the constrained PML estimates of $\beta$, say $\hat{\beta}^{r}_{T}$, and the Kuhn-Tucker (hereafter, KT) multipliers, say $\hat{\lambda}_{T}$, corresponding  to the inequality constraints. CFS demonstrate that one can use these linear combinations as intermediate statistics to produce I-I estimators of $\theta^0$ with asymptotically Gaussian limits. In addition, CFS show that imposing additional inequality restrictions on the auxiliary model will never decrease the efficiency of the resulting I-I estimator, so long as the corresponding KT multipliers are included in the vector of intermediate statistics. While correct, until now the reason behind this phenomena has not been completely understood. 

Our first contribution is to demonstrate that the asymptotic normality of the CFS I-I estimator is not due to the information brought by the constraints but is a direct consequence of the relationship between $\hat{\beta}^{r}_{T}$ and $\hat{\lambda}_{T}$. In particular, we demonstrate that the ``well-behaved'' linear combinations of $\hat{\beta}^{r}_{T}$ and $\hat{\lambda}_{T}$ put forward by CFS as new ``auxiliary
parameters" is asymptotically equivalent to a new feasible unconstrained PML estimator of the auxiliary model, which is always well-defined and asymptotically normal. Therefore, adding constraints to an auxiliary model does not increase the information about the structural parameters because this information was already contained in the (unconstrained) auxiliary model.\footnote{However, we must acknowledge that there are cases where equality constraints on the auxiliary parameters are indeed necessary to obtain a well-defined inverse of the binding function (see Gospodinov, Komunjer and Ng, 2017 for a recent  example).} 

Using this new feasible unconstrained estimator, we propose a computationally simple unconstrained I-I estimation strategy that does not enforce the inequality constraints on the auxiliary model. Even though our new I-I estimators do not enforce the inequality constraints on the auxiliary parameters, our I-I estimators are asymptotically equivalent to the constrained I-I estimators proposed by CFS. Moreover, we demonstrate that the standard asymptotic Gaussian distribution of our I-I estimators
remains valid, even if the {pseudo-true  value of the auxiliary parameters, $\beta^0$,} is on or near the boundary of the parameter space.

Our second contribution is to make rigorous the notion of pseudo-true values of the auxiliary parameters on or near the boundary of the auxiliary parameter space. While the approach of CFS treats the case of $\beta^0$ on the boundary, their approach is dichotomous: $\beta^0$ is either on the boundary, in which case $\hat{\beta}^r_{T}$ is not asymptotically normal, or it is not on the boundary, in which case $\hat{\beta}^r_{T}$ is asymptotically normal. This binary treatment of parameters near the boundary cannot capture cases where $\beta^0$ is close enough to the boundary of the parameter space to render standard asymptotic approximations unreliable. In such cases, modeling the pseudo-true value as a sequence that is close to but not on the boundary, for any finite sample size, can provide more reliable asymptotic approximations.

Following, among others, Moon and
Schorfheide (2009), {Andrews and Cheng (2012)}, and Ketz (2018), we consider a drifting
sequence of true data generating processes (DGPs), which, in turn, admits a sequence of drifting pseudo-true values for the auxiliary parameters. This sequence of pseudo-true auxiliary parameters is then employed to make rigorous the notion of parameters on or near the boundary within I-I estimation. This treatment requires us to revisit the limit theory of Andrews (1999) to accommodate cases where the drifting pseudo-true values cause a ``boundary bias'' that results in the score of the auxiliary model losing its asymptotic mean-zero property; we refer the reader to Section \ref{sec2} for precise details and discussion. However, even in the case where boundary bias occurs, our proposed I-I estimator still displays standard Gaussian asymptotics. The intuition behind this result is simple: I-I, by generating simulated data that mimics the observed data, carries out an implicit bias correction that alleviates the impact of the auxiliary parameters being close to the boundary.

We apply this new I-I approach to a range of examples that have featured in the I-I literature: a stochastic volatility model with a GARCH(1,1) auxiliary model (see, e.g., CFS); $\alpha$-stable models with a skewed Student-t auxiliary model (see, e.g., Garcia et al., 2011); and continuous-time jump-diffusion models for returns with a Student-t GARCH auxiliary model. In each example, we require that the auxiliary parameters satisfy a vector of inequality constraints, and, in each example, we demonstrate that empirically plausible values of the structural parameters lead to estimates of the auxiliary parameters that are near the boundary of the parameter space. We then demonstrate that our I-I approach can easily be applied to obtain estimators of the structural parameters that have good finite-sample performance.

The remainder of the paper is organized as follows. In Section two, we discuss constrained auxiliary models and give three classes of empirically relevant examples from the I-I literature where the constraints imposed on the auxiliary model are known to bind, at least in some cases. In addition, we present the particular drifting DGP framework considered in this paper, which generalizes the approach of Andrews (1999) to consider drifting pseudo-true values that can capture boundary affects for any finite-sample size. Within this particular setup, we demonstrate that a
well-defined unconstrained auxiliary parameter estimator, which contains the same
amount of information as the linear combinations of constrained auxiliary
estimates and KT multipliers used in CFS as auxiliary parameters, always exists and can be readily used for the purpose of I-I.
Section three uses this unconstrained auxiliary estimator to propose novel
I-I estimators and demonstrates that this
unconstrained estimator is asymptotically equivalent to the constrained I-I
approach proposed in CFS. In Section four we consider A series of Monte Carlo examples and an empirical application that demonstrates the performance of this approach and makes clear the empirical relevance of our approach. Section
five concludes and all proofs are relegated to the appendix.

\section{Inequality Constraints on the Auxiliary Model}\label{sec2}
We observe a sample $\{y_{T}:T\geq1\}$ generated from a strictly
stationary and ergodic probability model $P_{\theta }$ depending on an
unknown parameter $\theta \in \Theta\subset\mathbb{R}^{d_{\theta}}$, with $\Theta$ compact. Conditional on observed data $\mathbf{Y}_{t-1}=\{y_{t-1},y_{t-2},...\}$, the model admits the conditional density $%
p(y_{t}|\mathbf{Y}_{t-1};\theta )$. We are interested in conducting inference on $\theta$ in situations where maximum likelihood
estimation based on $p(y_{t}|\mathbf{Y}_{t-1};\theta )$ is infeasible or
otherwise unattractive, but simulation from $p(y_{t}|\mathbf{Y}_{t-1};\theta
)$ is relatively simple.

I-I proposes to estimate $\theta$ by targeting
consistent parameter estimates of a simpler auxiliary model $f(y_{t}|\mathbf{Y}_{t-1};\beta)$, with auxiliary parameters $\beta\in\mathbf{B}\subset\mathbb{R}^{d_{\beta}}$, with $\mathbf{B}$ compact, and where $d_{\beta}\geq
d_{\theta}$. Let $Q_{T}(\beta)$ {denote the sample auxiliary objective function
associated with $f(y_{t}|\mathbf{Y}_{t-1};\beta)$ and the observed sample $\{y_{t}\}_{t=1}^{T}$.}

We concern ourselves with situations where, to ensure estimates obtained from $Q_{T}(\beta)$ are well-behaved, $\beta$ must be estimated subject to a vector of inequality restrictions: $$g(\beta)\geq 0.$${The constraint function $g:\breve{\mathbf{B}}\rightarrow 
\mathbb{R}^{q}$, with $\breve{\mathbf{B}}\subset\mathbb{R}^{d_\beta}$ an open set containing $\mathbf{B}$,} is known and continuously
differentiable on $\breve{\mathbf{B}}$. Throughout the remainder the notation $g(\beta )\geq 0$ {is taken to mean} $g_{j}(\beta)\geq 0$, $j=1,...,q$. From the inequality constraints, we define a restricted, or constrained, parameter space as $$\mathbf{B}^r:=\{\beta\in\mathbf{B}:g(\beta)\geq 0\}.$$ {Throughout, we allow the number of constraints $(q)$ to be greater than or less than the number of parameters $(d_\beta)$. }

We assume throughout that $\mathbf{B}^r$ has a non-empty interior, which precludes equality constraints. However, this assumption is immaterial since, up to an abuse of notation, if the problem originally featured a mix of equality constraints, say $\{g_l(\beta)=0:1\leq l\leq q_{1}\}$, and inequality constraints, say $\{g_k(\beta)\geq 0:1\leq k\leq q_{2}\}$, with $q_{1}+q_{2}= q$,  we could always re-define $g(\beta)\geq 0$ to be only the inequality constraints $\{g_k(\beta)\geq 0:1\leq k\leq q_{2}\}$ that remain active after imposing the equality constraints and eliminating some corresponding components of $\beta$.

Subsequently, we can define the constrained estimator of the auxiliary parameters as 
\begin{flalign*}
\hat{\beta}^r_T&:=\arg\max_{\beta\in\mathbf{B}}Q_{T}(\beta)\text{ s.t. }g(\beta)\geq 0.\\
&:=\arg\max_{\beta\in\mathbf{B}^r}Q_{T}(\beta)
\end{flalign*}
Consistent I-I estimation requires that $\hat{\beta}^r_{T}$ be a consistent estimator of an appropriately defined pseudo-true parameter value, generically denoted by $\beta ^{0}$ and satisfying $g(\beta^0)\geq 0$.  As discussed by CFS, if these inequality constraints bind at $\beta^0$, in the sense that for at least one $j$,  $$g_j(\beta^0)=0,\;\;1\leq j\leq q,$$I-I based on $\hat{\beta}^r_T$ may result in estimators with a non-Gaussian limit distribution.

In certain examples, such as those considered in the following section, in order for estimators based on $Q_{T}(\beta)$ to be well-behaved, {certain inequalities must be strict, in that $g_j(\beta)>0$ for some $1\leq j\leq q$.} We will ensure that $\hat{\beta}^r_{T}$ satisfies this property, at least with probability one for $T$ large enough, by assuming that we have a drifting DGP where, for any finite $T$, the pseudo-true value of the optimization program belongs to $\text{Int}(\mathbf{B}^r)$. In so doing, when some constraints are actually binding, this precisely means that the pseudo-true value is then ``near the boundary''; we refer to, e.g., Andrews and Cheng (2012), and Ketz (2018) for a similar use of this terminology.  

Before going further with the precise mathematical framework, we first give examples where constraints on the auxiliary models used within I-I estimation feature in the empirical literature.

\subsubsection*{Example 1: Stochastic Volatility}\label{examp}

We begin with the classic example of a log-normal
stochastic volatility (SV) model and a GARCH(1,1) auxiliary model with Gaussian or Student-t innovations.\footnote{Several authors consider I-I estimation for this model using GARCH$(p,q)$ models and we refer the reader to Engle and Lee (1996), Monfardini (1998),
	Pastorello et al. (2000) and CFS for examples.} The log-normal stochastic volatility model is defined as follows:
\begin{eqnarray}
y_{t} &=&\sqrt{h_{t}}e_{t},\;t=1,...,T  \label{SV}, \\
\ln (h_{t}) &=&\alpha +\delta \ln (h_{t-1})+\sigma _{v}v_{t},  \nonumber
\end{eqnarray}
where $|\delta |<1$, $\sigma _{v}>0$, $(e_{t},v_{t})^{\prime }\sim
_{i.i.d.}N(0,\text{Id}_{2})$ and we denote the structural parameters as $%
\theta =(\alpha ,\delta ,\sigma _{v})^{\prime }$. We observe a series $%
\{y_{t}\}_{t=1}^{T}$ from the SV model in (\ref{SV}) and our goal is to
conduct inference on $\theta $. 

We follow CFS and consider as our auxiliary model the
GARCH(1,1) model:
\begin{eqnarray}
y_{t} &=&\sqrt{h_{t}}\epsilon_{t}  \label{GARCH}, \\
h_{t} &=&\psi +\varphi y_{t-1}^{2}+\pi h_{t-1}.  \nonumber
\end{eqnarray}
Common specifications for the errors $\epsilon _{t}$ in (\ref{GARCH}) are $%
\epsilon _{t}\sim_{iid} N(0,1)$ {or, for $v(\eta):=[{{(1/\eta-1/2)}/{1/\eta}}]$, $\epsilon _{t}\sim_{iid} v(\eta)^{1/2}t_{1/\eta}$, where $t_{1/\eta}$ denotes the Student-t distribution with $1/\eta $ degrees of freedom, so that $v(\eta)^{1/2}t_{1/\eta}$ denotes a Student-t with unit variance (which requires that $\eta<1/2$). }We denote the
auxiliary parameters as $\beta $, with $\beta =(\psi
,\varphi ,\pi )^{\prime }$ if $\epsilon _{t}\sim_{iid} N(0,1)$
and $\beta =(\psi ,\varphi ,\pi ,\eta )^{\prime }$ otherwise. The GARCH(1,1)
model is very useful as an auxiliary model as it can capture many of the
structural ideas associated with (\ref{SV}), such as thick tails and
volatility clustering, while yielding closed form formulas for the score and
Hessian based on the pseudo-log-likelihood $Q_{T}(\beta )$. 

However, the GARCH(1,1) auxiliary model must be estimated subject to inequality
constraints to ensure that the pseudo-maximum likelihood estimator of $\beta$ is well-behaved. The set of inequality constraints for the auxiliary model can be stated
as
\begin{equation}
\psi \geq 0,\;\varphi \geq0,\;\pi \geq 0,\;\varphi +\pi \leq 1,
\label{GARCHineq}
\end{equation}%
with the added constraint $0\leq \eta \leq1/2$ when $\epsilon _{t}$ is distributed as Student-t with $1/\eta $ degrees of freedom. Moreover, we stress that the quasi-likelihood is not even well-defined when all the parameters, $(\varphi,\psi,\pi)'$ are simultaneously on the boundary, since the conditional variance must be strictly positive. More generally, Francq and Zakoian (2007) stress that to establish asymptotic normality of the QMLE for the GARCH parameters, a key regularity condition is that ``the true parameter must lie in the interior of the parameter space''. This statement can obviously be extended to the pseudo-true value of the GARCH parameters. To
enforce the above strict inequalities on the auxiliary parameters, CFS require
(see their footnote five on page 960) that the GARCH parameters in their auxiliary
model satisfy
\begin{equation}
\varphi \geq 0.025,\;\eta \leq 0.499 . \label{CFStrick}
\end{equation}
However, such constraints are effectively arbitrary and do not necessarily reflect the true nature of the constraints in \eqref{GARCHineq}.


\subsubsection*{Example 2: $\alpha$-Stable Random Variables}
The class of $\alpha$-stable distributions is often used to capture random variables that display heavy-tailed features, such as stock returns data. The distribution of a random variable, $y_t$, from the $\alpha$-stable class is characterized by four parameters: $\alpha$- the tail index, which captures the `heavyness' of the tail; $\gamma$- the skewness parameter; $\mu$- the location parameter; and $\sigma$- the scale parameter. Denoting $\theta=(\alpha,\gamma,\mu,\sigma)'$, we have that $$\Theta:=(0,2)\times(-1,1)\times\mathbb{R}\times[0,\infty).$$ {We note here that values of $\alpha<2$ ensure that the variance of the random variable is not finite, while if $\alpha\leq 1$ both the mean and the variance are not finite.} We refer the reader to Samoradnitsky (2017) for a book length treatment on $\alpha$-stable random variables. 

An interesting feature of the $\alpha$-stable class is its lack of a closed-form density function, which makes the application of maximum likelihood methods to estimate $\theta$ difficult (see, e.g., {Garcia et al., 2011} for a discussion). The difficulty of maximum likelihood estimation has led to the development of I-I estimators for $\theta$ that first postulate an auxiliary model with parameters $\beta$ that can roughly match the parameters of the $\alpha$-stable distribution and for which consistent estimators of these parameters can easily be obtained.

Following Garcia et al. (2011), one such class of auxiliary models is the skewed Student-t (hereafter, skew-t) distribution developed by {Fernandez and Steel (1998)}: 
\begin{flalign*}
f(y;\beta)=\frac{\frac{\Gamma(\frac{\nu}{2}+\frac{1}{2})}{\Gamma(\nu/2)}}{\sqrt{\nu}}\frac{1}{\ell\left(\eta+\frac{1}{\eta}\right)}\left\{1+\frac{1}{\nu}\left(\frac{y-\omega}{\ell}\right)^2\left[\frac{1}{\eta^2}\1[y\geq \omega]+\eta^2\1[y<\omega]\right]\right\}^{-\frac{\nu+1}{2}}.
\end{flalign*}The degree of freedom parameter $\nu$ captures tail thickness, $\eta$ captures skewness, and $\omega$ and $\ell$ denote the location and scale parameters. Clearly, the parameters $\theta$ and $\beta$ are closely related. Moreover, the close match between the parameters of the $\alpha$-stable and skew-t distributions should lead to well-behaved and nearly efficient I-I estimators. 

{In many empirical applications, estimates of $\alpha$ are often near 2; recall that a value of $\alpha<2$ implies that the unconditional variance is not finite.} This feature is potentially troubling for I-I estimation since Garcia et al. (2011) demonstrate that when $\alpha$ is larger than about 1.9, the parameter $\nu$ in the skew-t auxiliary model, {and by association the $\alpha$ parameter in the structural model, becomes poorly identified. In particular, Garcia et al (2011) argue that ``for $\alpha$ close to 2, we may expect that observed data will give the spurious
feeling that variance is finite, which would imply a normal distribution corresponding to
$\nu = +\infty$ in a Student framework. This is why we will constrain the auxiliary parameter.'' Since the authors assume that $\alpha<2$, the authors constrain the skew-t auxiliary model in a similar fashion and impose the inequality constraint $\nu\leq 2$ on the auxiliary parameter.} Numerical results presented in Garcia et al. (2011) demonstrate that a constrained version of I-I, which uses the inequality constraint $\nu\leq2$, produces estimators that are better behaved than those where the constraint on the auxiliary parameter $\nu$ is not maintained. A similar I-I estimation strategy has also been employed by {Lombardi and Calzolari (2009)} for estimation of $\alpha$-stable stochastic volatility models and by {Calzolari and Halbeib (2018)} for estimation of $\alpha$-stable factor models.

\subsubsection*{Example 3: Stochastic Volatility Jump-Diffusion (SVJD) Models}
The stylized facts of time-varying and autocorrelated volatility, allied with non-Gaussian return distributions, are now extensively documented in the literature on financial returns. However, it is often the case that standard volatility models, such as the one treated in Example 1, cannot completely capture the variability of daily returns in periods of extreme volatility, such as during the 2008-2009 financial crisis. A common approach to address this issue is to consider the inclusion of so-called `jump' processes within existing volatility models. The inclusion of the jump process allows volatility models to exhibit periods of high volatility without significantly altering the interpretation of the model. The literature on modeling returns under the assumption of non-negligible jumps is now extensive, and we refer the reader to {Ait-Sahalia and Jacod (2014)} for a textbook treatment.

 An important class of widely used continuous-time stochastic volatility models in finance is the mean reverting stochastic volatility jump diffusion (SVJD) model.
Let $P_{t}$ denote the asset price at time $t>0$, and let $p_{t}:=\ln(P_{t})$. In the SVJD model, the evolution of $p_t$ follows a bivariate jump diffusion process, with a representative example of the SVJD model being
\begin{flalign} d p _ { t } & =  { \mu } d t + \exp ({ V _ { t }/2 }) d W _ { t } ^ { p } + d J _ { t } ^ {  }, \nonumber\\ d V _ { t } & =  { \kappa } \left(  { \eta } - V _ { t } \right) d t +  { \sigma } _ { v }  d W _ { t } ^ { v }, \label{svj}
\end{flalign}where $dW_{t}^{v},\;dW_{t}^{p}$ are independent standard Brownian motion processes, $d J_{t}$ is a jump component with $d J_{t}:=Z_{t}dN_{t}$, $Z_{t}$ denotes the jump size and $dN_{t}$ is a Poisson process with constant jump intensity. The SVJD model captures two important empirical features of asset prices: one, return volatility exhibits strong serial dependence; two, price jumps exist due to the arrival of unanticipated market news. 

The SVJD model is a workhorse of empirical finance, where nonparametric approaches are commonly used to obtain high-frequency measures of variability, such as integrated volatility and quadratic variation (we refer the reader to {Andersen et al., 2009} for a discussion of volatility measures in this model, and its many generalizations).  On the other hand, inference on the unknown parameters is hindered by the latent nature of the volatilities, which ensures that estimation techniques based on the likelihood are computationally demanding. Luckily, simulation-based procedures, such as I-I, can bypass the calculation of the likelihood by simulating directly from the process in \eqref{svj}{; for empirical applications of the SVJD model using simulation-based inference techniques see, e.g., Eraker (2001), Andersen et al. (2002), Creel and Kristensen (2015).} 
 

For the purpose of I-I estimation, a useful class of auxiliary models for capturing the behavior in \eqref{svj} would be the GARCH model class in equation \eqref{GARCH}. As discussed in Example 1, the use of GARCH auxiliary models requires imposing several inequality restrictions. However, even for this simple version of the SVJD model in equation \eqref{svj}, GARCH auxiliary models can often yield estimates of the auxiliary parameters for which the inequality constraints in \eqref{GARCHineq} bind. For example, GARCH-based estimators of $\varphi+\pi$ are often very close to unity.\footnote{We refer the reader to the Monte Carlo section for numerical evidence of this statement as it pertains to the GARCH auxiliary model in equation \eqref{GARCH}.} Therefore, as argued in CFS, additional information about the auxiliary parameters would be required to successfully identify the structural parameters.

\subsection{Parameters Near the Boundary}Before presenting our new approach to I-I with constraints, we clarify what is meant by auxiliary parameters on the boundary of the parameter space. Recall the definition $$\mathbf{B}^r:=\{\beta\in\mathbf{B}:g(\beta)\geq 0\}.$$ We say that auxiliary parameters are on the {boundary if $\beta^0$, the pseudo-true value of the auxiliary parameters, is on the boundary of $\mathbf{B}^r$. }

{We rely on a drifting DGP to capture the behavior of
extremum estimators when $\beta^0$ is on the boundary of $\mathbf{B}^r$.} In particular, we consider that the DGP of the structural model is indexed by a sequence of drifting ``true'' values $\{\theta_{T}\}:=\{\theta_{T}:T\geq1\}$ that satisfy\footnote{{Formally, this assumption implies that the observed data, viewed as a triangular array $\{y_{t,T}:T\geq1,t=1,\dots,T\}$, comes from a stationary process depending on $T$.}} $$\{\theta_T\}\in{\Theta}[\mathbb{N}^{+}]:=\left\{\{\theta_T\in\Theta:T\geq1\}:\;\lim_{T\rightarrow\infty}\theta_T=\theta^0\in\text{Int}(\Theta)\right\}.$$  The population objective function for the auxiliary model, calculated under this DGP, is denoted by $\mathcal{Q}(\theta_T,\beta)$. Using $\mathcal{Q}(\theta_T,\beta)$, we define the sequence of pseudo-true auxiliary parameters $$b(\theta_{T}):=\arg\max_{\beta\in\mathbf{B}^r}\mathcal{Q}(\theta_T,\beta),\text{ where }\{\theta_T\}\in{\Theta}[\mathbb{N}^{+}].$$ The notation $\theta_T\mapsto b(\theta_T)$ clarifies that this map depends on the sequence $\{\theta_T\}$. To maintain notational simplicity, when no confusion is likely to result, we denote these pseudo-true auxiliary parameters as $\beta^0_T:=b(\theta_T)$. 

The drifting sequences $\{\theta_T\}$ and $\{\beta^0_T\}$ allow us to  capture auxiliary parameters near the boundary of $\mathbf{B}^r$ using the set: 
\begin{flalign*}
\Gamma(\theta^0,\beta^0)&:=\left\{\{\theta_T\}\in {\Theta}[\mathbb{N}^{+}]:\beta^0_T\in\text{Int}(\mathbf{B}^r),\;\lim_{T\rightarrow\infty}\beta^0_T:=\beta^0\in \mathbf{B}^{r}\right\},
\end{flalign*} and by restricting our analysis to DGPs satisfying \begin{flalign}\label{driftDGP} 
\{\theta_T\}\in\Gamma(\theta^0,\beta^0).
\end{flalign}
This construction enforces that
$\beta^0 _{T}$ belongs to the interior of $\mathbf{B}^r$ but allows $\beta^0$ to lie on the boundary of $\mathbf{B}^r$.\footnote{We note that, when $\beta^{0} $\ is in the interior of the parameter set, the concept of a drifting DGP is
	hardly useful; one can then assume $b(\theta _{T})=\beta ^{0}$\ for all $T$\
	sufficiently large. For instance, in the illustrative stochastic volatility example, a drifting true value for $\psi^{0}$ is not necessary; when $\varphi^{0}$ is on the boundary ($\varphi^{0}=0$), $\psi^{0}$ must be strictly positive and its constrained estimator (the sample mean of $y_{t}^{2}$) will automatically fulfill this inequality constraint. } 

Our need to consider a drifting DGP partly arises because we have in mind cases when $\beta $ must satisfy some strict inequalities in order for $Q_{T}(\beta )$ to be
well-behaved in finite-samples, but where the population analogue $\mathcal{Q}(\theta_T,\beta)$ remains well-behaved for all $\beta\in \mathbf{B}^r$. The definition of $\Gamma(\theta^0,\beta^0)$ ensures this by requiring that $\beta _{T}^{0}\in\text{Int}(\mathbf{B}^r)$, so that $Q_{T}(\beta^0_{T} )$ remains well-defined, while $\lim_T\mathcal{Q}(\theta_T,\beta^0_T)$ is well-defined under standard continuity assumptions. As an illustration of why we require this drifting DGP, consider the stochastic volatility example, and recall that the auxiliary parameter $\varphi$ {must be strictly positive to ensure that the pseudo-true value of $\beta$ is identified}. To enforce this condition CFS require the ad hoc condition \eqref{CFStrick}, while we enforce this condition by
imposing the high-level condition \eqref{driftDGP}.

It is worth noting that we maintain the assumption that $b(\theta _{T})$\ always fulfills the constraints, and, by
continuity, $\beta ^{0}$ must also fulfill them, while it may violate the strict inequality constraints we implicitly need to
maintain in finite samples. However, we may expect that
all KT multipliers still converge to zero, in contrast to the
setting considered {in the asymptotic theory of CFS}. 


\subsection{Standard Auxiliary Estimators}
I-I with constrained auxiliary parameters relies on the constrained estimator $\hat{\beta}^{r}_{T}$ which satisfies
\[
Q_{T}(\hat{\beta}_{T}^{r})\geq \sup_{\beta \in \mathbf{B}^r}Q_{T}(\beta )+o_{P}(1/T),
\]%
which we can obtain through the Lagrangian function
\[
\mathcal{L}_{T}(\beta ,\lambda ):=Q_{T}(\beta )+g(\beta )^{\prime }\lambda ,
\]for $\lambda\in\mathbb{R}^q$ a vector of KT multipliers. Under differentiability
conditions, $\hat{\beta}_{T}^{r}$ and the associated KT
multipliers $\hat{\lambda}_{T}$\ solve the first-order
conditions
\begin{equation}
\frac{\partial Q_{T}(\hat{\beta}_{T}^{r})}{\partial \beta }+\frac{\partial
	g^{\prime }(\hat{\beta}_{T}^{r})}{\partial \beta }\cdot \hat{\lambda}_{T}=0
\label{lagrange},
\end{equation}
with the slackness conditions
\begin{flalign}
\hat{\lambda}_{j,T}\cdot g_{j}(\hat{\beta}_{T}^{r})&=0,\text{ for all } j=1,...,q
\label{slackness} \\
g(\hat{\beta}_{T}^{r}) &\geq 0,\;\hat{\lambda}_{T}\geq 0. \nonumber
\end{flalign}

Building on the theory of constrained estimation discussed in Andrews (1999), under the following sufficient conditions, the estimators $\hat{\beta}^r_{T}$ and $\hat{\lambda}_{T}$ are $\sqrt{T}$-consistent estimators of $\beta^0_{T}$ and $0$, respectively.\footnote{These assumptions are similar to those given in CFS (see their Assumptions 1 and 3), but are adapted to accommodate our drifting DGP setting.}
	
	\medskip 
	
		\noindent\textbf{Assumption A0:} Under $\{\theta_T\}\in\Gamma(\theta^0,\beta^0)$: (i) $\sup_{\beta\in\mathbf{B}^r}|Q_T(\beta)-\mathcal{Q}(\theta_T,\beta)|=o_{P}(1);$ (ii) For all $\varepsilon>0$ , $\liminf_{T}\{\mathcal{Q}(\theta_T,\beta^0_T)-\sup_{\beta\in\mathbf{B}^r:\|\beta-\beta^0_T\|>\varepsilon}\mathcal{Q}(\theta_T,\beta)\}>0$.

\medskip

\noindent\textbf{Assumption A1:} Under $\{\theta_T\}\in\Gamma(\theta^0,\beta^0)$:

\noindent (i) $\beta \mapsto Q_{T}(\beta )$ has continuous partial derivatives of order two on $\text{Int}(\mathbf{B}_{}^{r})$ with probability one. 

\medskip 

\noindent (ii) For $\mathcal{J}^{0}$ a non-stochastic $(d_{\beta }\times d_{\beta
}) $\ positive-definite matrix, and for any $\gamma>0$: $$\sup_{\beta\in\mathbf{B}^{r}:\|\beta-\beta^0_{T}\|\leq\frac{\gamma}{\sqrt{T}}}\left\|\frac{\partial^2 Q_{T}(\beta)}{\partial\beta\partial\beta'}+\mathcal{J}^0\right\|=o_{P}(1). $$

\noindent(iii) For some $\delta^0\in\mathbb{R}^{d_{\beta}}$, and for $\mathcal{I}^0$ a non-stochastic $(d_{\beta }\times d_{\beta
}) $\ positive-definite matrix, $$\sqrt{T}{\partial Q_{T}(\beta
	_{T}^{0})}/{\partial \beta ^{ }}\rightarrow _{d}\aleph \left(\delta^0,\mathcal{I}^{0}\right).$$

\noindent(iv) {(a) ${g}(\beta)$ is continuously differentiable for $\beta\in\breve{\mathbf{B}}$, for $\breve{\mathbf{B}}$ open and $\mathbf{B}\subseteq\breve{\mathbf{B}}$; (b) there exists a $0\leq\tilde{q}\le \min\{q,d_\beta\}$ and a function $\tilde{g}:\breve{\mathbf{B}}\rightarrow\mathbb{R}^{\tilde{q}}$, a $\tilde{q}$-dimensional sub-vector of $g(\beta)$, such that, for all $T$ large enough, $\tilde{g}(\beta^0_T)$ contains all the zero entries of $g(\beta^0_T)$ and  $\text{rank}\left(\partial \tilde{g}(\beta^0_T)'/\partial\beta\right)=\tilde{q}$.}
\medskip

\medskip 

\noindent\textbf{Assumption A2:} Consider the quadratic expansion: 
	\begin{equation*}
	Q_{T}(\beta )=Q_{T}(\beta _{T}^{0})+\frac{\partial Q_{T}(\beta _{T}^{0})}{%
		\partial \beta ^{\prime }}(\beta -\beta _{T}^{0})+\frac{1}{2}(\beta -\beta
	_{T}^{0})^{\prime }\frac{\partial ^{2}Q_{T}(\beta _{T}^{0})}{\partial \beta
		\partial \beta ^{\prime }}(\beta -\beta _{T}^{0})+R_{T}(\beta ) . 
	\end{equation*}
Under $\{\theta_T\}\in\Gamma(\theta^0,\beta^0)$: for any sequence $\gamma _{T}=o(1)$,
\[
\sup_{\beta\in\mathbf{B}^{r}:\left\Vert \beta -\beta _{T}^{0}\right\Vert
	\leq \gamma _{T}} \left\{ \frac{\left\vert  R_{T}(\beta )\right\vert }{\left[
	1+\sqrt{T}\left\Vert \beta -\beta _{T}^{0}\right\Vert \right] ^{2}}\right\}
=o_{P}\left(1/T\right) .
\]

\medskip

The above assumptions are similar to those employed by Andrews (1999) to deduce his Theorem 1. However, Assumption \textbf{A1}(iii) is novel and is maintained to accommodate our drifting DGP. In particular, 
Assumption \textbf{A1}(iii) allows for the drifting behavior of  $\beta^0_T$ to contaminate the limiting distribution of the scaled pseudo-score $\sqrt{T}\partial Q_{T}(\beta^0_T)/\partial\beta$. In particular, Assumption \textbf{A1}(iii) can capture cases where the proximity of $\beta^0_T$ to the boundary causes a ``boundary bias'', whereby the pseudo-score looses its asymptotic mean-zero property, typically because $\sqrt{T}(\beta^0_T-\beta^0)$ {is $O(1)$ and not $o(1)$.} This framework is often used when one wishes to accurately capture the behavior of estimators when the estimated parameter values are close to, but potentially not on, the boundary of the parameter space; while we have suggested several examples where this phenomena may be in evidence, we refer the interested reader to Ketz (2018) for further examples and discussion. {Assumption \textbf{A1}(iv) deals with the behavior of the constraint function $g(\beta)$. Part (a) is standard, while part (b) is required since we allow the number of inequality constraints to be larger than the number of auxiliary parameters (i.e., $q>d_\beta$). Indeed, one particularly important example of this phenomena is the GARCH auxiliary model discussed in Section two. Intuitively, $\tilde{g}(\beta)$ includes all the constraints that are active at $\beta^0_T$.}

Under the above assumptions, the following result holds. 


\medskip 

\noindent\textbf{Lemma 1:} For $\{\theta_T\}\in\Gamma(\theta^0,\beta^0)$, under \textbf{A0-A2}, $\sqrt{T}(\hat{\beta}^r_T-\beta^0_T)=O_{P}(1)$ and $\sqrt{T}\hat{\lambda}_{T}=O_{P}(1)$. 

\medskip

While $\hat{\beta}^{r}_{T}$ is $\sqrt{T}$-consistent, it has rightly been stressed in CFS that the constraints $g(\beta)\geq 0$ may produce some singularity
(and non-normality) in the asymptotic distribution of $\hat{\beta}_{T}^{r}$, and therefore $\hat{\beta}^r_{T}$ may not be appropriate for I-I estimation. {This issue is exacerbated under our setup, since not only may the asymptotic distribution of $\sqrt{T}(\hat{\beta}^r_{T}-\beta^0_T)$ be non-normal, but, under Assumption \textbf{A1}(iii), the asymptotic distribution of $\sqrt{T}(\hat{\beta}^r_{T}-\beta^0_T)$ may have a non-zero asymptotic mean.} For the former reason, i.e., non-normality, CFS search for a seemingly ad
hoc linear combination of the constrained estimator $\hat{\beta}_{T}^{r}$
and the vector $\hat{\lambda}_{T}$ of KT multipliers that is asymptotically
normal (see Proposition 2 in CFS, page 950). {To elucidate the implications of this statement, first consider the following (infeasible) quadratic objective function based on the quadratic expansion in Assumption \textbf{A2}:}
\[
{M}_{T}(\beta )=Q_{T}(\beta _{T}^{0})+\frac{\partial Q_{T}(\beta
	_{T}^{0})}{\partial \beta ^{\prime }}\left(\beta -\beta _{T}^{0}\right) +%
\frac{1}{2}\left(\beta -\beta _{T}^{0}\right)^{\prime }\frac{\partial
	^{2}Q_{T}\left( \beta _{T}^{0}\right) }{\partial \beta \partial \beta
	^{\prime }}\left( \beta -\beta _{T}^{0}\right).
\]{Since $\beta _{T}^{0}$\ is in the
interior of the parameter space, $M_T(\beta)$ is
well-defined for any $\beta \in 
\mathbb{R}
^{d_{\beta }}$
and is uniquely maximized by the solution of the first-order conditions,
}
\[
\ddot{\beta}_{T}=\beta _{T}^{0}+J_{T}^{-1}\frac{\partial Q_{T}(\beta
_{T}^{0})}{\partial \beta },\text{ where }J_{T}=-\frac{\partial
^{2}Q_{T}\left( \beta _{T}^{0}\right) }{\partial \beta \partial \beta
^{\prime }}.
\]

Note that, with an abuse of language, $\ddot{\beta}_{T}$\ could be dubbed an
``unconstrained estimator" since the constraints $g\left( \beta \right) \geq 0
$\ are never taken into account in its definition. The advantage of this
``estimator" is that it always exists, since $\beta^0_T$ is an interior point the quadratic approximation always exists, while a general unconstrained estimator may not
even exist. However, calling $\ddot{\beta}_{T}$\ an estimator is an abuse of language since it is not feasible to
{compute $M_{T}(\beta)$ since the pseudo-true value $\beta^0_T$ is unknown.}

Our first key result of this section is to demonstrate that, under our drifting DGP setup, the
linear combination of auxiliary parameters put forward in CFS is tightly related to this
(potentially) infeasible unconstrained estimator. 

\medskip

\noindent \textbf{Proposition 1:} \noindent For $\{\theta_T\}\in\Gamma(\theta^0,\beta^0)$, under Assumptions \textbf{A0}-%
\textbf{A2}:%
\begin{equation}
J_{T}\sqrt{T}\left( \hat{\beta}_{T}^{r}-\beta _{T}^{0}\right) -\frac{%
\partial g^{\prime }(\beta _{T}^{0})}{\partial \beta }\sqrt{T}\hat{\lambda}%
_{T}=J_{T}\sqrt{T}\left( \ddot{\beta}_{T}-\beta _{T}^{0}\right) +o_{P}(1)
\label{combinCFS}.
\end{equation}
The remainder term $o_{P}(1)$ in (\ref{combinCFS}) is identically zero when
the criterion function $Q_{T}(\beta )$\ is quadratic and the constraints $%
g(\beta )$ are linear.\hfill $\square $ \medskip 

The LHS of equation (\ref{combinCFS}) is identical to the so-called
\textquotedblleft linear combinations [of the constrained estimator and KT
multipliers] that are asymptotically well behaved" in Proposition 2 of CFS
(pg 950). By \textquotedblleft well-behaved" CFS essentially mean
asymptotically normal, whereas separately the constrained estimator and the
KT multipliers may not be asymptotically normal\ when the parameters are
close to or on the boundary of $\mathbf{B}^{r}$.\footnote{When the constraints $g(\cdot )$ are non-linear, CFS actually
consider more complicated linear combinations involving the second
derivatives of $g(\beta )$. However, these additional
terms will cancel out when working, as we do in this section, under
the assumption that the constraints are fulfilled {in the population;} in this case, the
vector of KT multipliers actually converge to zero, and kill the
additional terms in CFS.}
Proposition 1 demonstrates that the linear combinations studied in CFS are
well-behaved, precisely because they correspond (asymptotically) to the
unconstrained extremum estimator; i.e., equation (\ref{combinCFS})
demonstrates that by combining the auxiliary parameters and KT multipliers
to create the well-behaved linear combinations, we are just back to
unconstrained estimation!

\subsection{Asymptotically Normal Feasible Unconstrained Estimation}
Proposition 1 demonstrates that the linear combinations of $\hat{\beta}_{T}^{r}$ and $\hat{\lambda}_{T}$ that lead to asymptotically normal
auxiliary parameters are asymptotically equivalent to the infeasible unconstrained
estimator $\ddot{\beta}_{T}$. Therefore, a feasible version of $\ddot{\beta}_{T}$ would provide an asymptotically equivalent alternative to the ad hoc
combination of constrained estimators and KT multipliers used in CFS. To
deduce such an estimator, we first recall that $\ddot{\beta}_{T}$ is
actually the global maximizer of {the quadratic 
objective function $M_{T}(\beta)$, which depends on the infeasible $\beta _{T}^{0}$}. {This suggests that a feasible unconstrained estimator can be obtained by replacing} $\beta _{T}^{0}$ in {$M_{T}(\beta)$ by a consistent estimator.} By Lemma 1, the constrained
estimator $\hat{\beta}_{T}^{r}$ is a consistent estimator of $\beta _{T}^{0}$%
, and we can replace $\beta _{T}^{0}$ in the quadratic objective function $M_T(\beta)$ by $\hat{%
\beta}_{T}^{r}$, and define the feasible unconstrained estimator 
\[
\widehat{\beta }_{T}=\arg \max_{\beta \in 
\mathbb{R}
^{d_{\beta }}}\left[ Q_{T}(\hat{\beta}_{T}^{r})+\frac{\partial Q_{T}(\hat{%
\beta}_{T}^{r})}{\partial \beta ^{\prime }}(\beta -\hat{\beta}_{T}^{r})+%
\frac{1}{2}(\beta -\hat{\beta}_{T}^{r})^{\prime }\frac{\partial ^{2}Q_{T}(%
\hat{\beta}_{T}^{r})}{\partial \beta \partial \beta ^{\prime }}(\beta -\hat{%
\beta}_{T}^{r})\right] .
\]%
Interestingly, $\widehat{\beta }_{T}$ is obtained simply by taking a
Newton-step away from $\hat{\beta}_{T}^{r}$: 
\[
\widehat{\beta }_{T}=\hat{\beta}_{T}^{r}-\left[ \frac{\partial ^{2}Q_{T}(%
\hat{\beta}_{T}^{r})}{\partial \beta \partial \beta ^{\prime }}\right] ^{-1}%
\frac{\partial Q_{T}(\hat{\beta}_{T}^{r})}{\partial \beta ^{{}}},
\]%
so that obtaining $\widehat{\beta }_{T}$ is extremely simple in practice.
Throughout the remainder, we refer to $\widehat{\beta }_{T}$ as the feasible
unconstrained (FUNC) estimator of $\beta _{T}^{0}$.

Before the FUNC estimator $\widehat{\beta }_{T}$ can be used for the purpose
of I-I, we must understand its asymptotic properties. The asymptotic
behavior of $\widehat{\beta }_{T}$ can be determined by analyzing the
asymptotic behavior of the quadratic expansion in Assumption \textbf{A2}. 	We now give the main result of this section: the FUNC estimator $\widehat{%
\beta }_{T}$ is asymptotically equivalent to $\ddot{\beta}_{T}$, and thus to
the well-behaved linear combinations employed in CFS as auxiliary parameters.

\medskip

\noindent \textbf{Theorem 1:} For $\{\theta_T\}\in\Gamma(\theta^0,\beta^0)$, under Assumptions \textbf{A0}-\textbf{A2},  $ 
\sqrt{T}\left( \widehat{\beta }_{T}-\ddot{\beta}_{T}\right) =o_{P}(1). 
\hfill\square$

	\medskip

		Ketz (2018) has proven a similar result in the framework of a drifting
	true value similar to ours. {For the sake of being self-contained, we provide our own proof of this result.	Before concluding, we note that, by the result of Theorem 1, the FUNC estimator $\widehat{%
		\beta }_{T}$\ allows us to
	rewrite the decomposition (\ref{combinCFS}) as follows:%
	\begin{equation}\label{13++}
	J_{T}\sqrt{T}\left( \hat{\beta}_{T}^{r}-\beta _{T}^{0}\right) -\frac{%
		\partial g^{\prime }(\beta _{T}^{0})}{\partial \beta }\sqrt{T}\hat{\lambda}%
	_{T}=J_{T}\sqrt{T}\left( \widehat{\beta }_{T}-\beta _{T}^{0}\right)
	+o_{P}(1) .
	\end{equation}
	That is, by working with the computationally friendly FUNC estimator $\widehat{\beta }%
	_{T}$ we convey exactly the same information as the complicated linear combination of constrained estimators and
	KT multipliers considered by CFS. The implications of this remark
	for the purpose of I-I are discussed in the subsequent
	sections.

\section{Indirect Inference With(Out) Constraints}\label{sec3}

The key input of I-I is a\ set of $H$\ simulated paths $\{%
\tilde{y}_{t}^{(h)}(\theta )\}_{t=1}^{T},h=1,..,H$. From this input, there
are several ways to perform I-I. Our focus of interest in this section is to compare
four strategies. The first two strategies are based on the score matching approach
of GT. The approach of CFS and the approach proposed
in this paper will produce two distinct, albeit asymptotically equivalent,
variants of the score-matching approach. As already mentioned in the
comments of Theorem 1, we differ from CFS in that we will not incorporate, explicitly, the KT multipliers as additional
auxiliary parameters for I-I since the FUNC estimator $\widehat{\beta}_{T}$\ carries
the same information.

The last two strategies are based on the GMR approach of minimum distance
between auxiliary parameters. These two strategies differ regarding the
parameters to match: constrained estimators of $\beta $\ augmented by KT multipliers, as in CFS, or the user-friendly FUNC estimator proposed
in this paper.

By analogy with the trinity of tests, we will dub ``Wald approach" the
minimum distance approach while the score-matching approach will simply be
called ``Score approach." Note that CFS dub CMD (Classical Minimum Distance)
the Wald approach and GMM (Generalized Method of Moments) the Score
approach.\ GMR have shown that in classical circumstances (I-I without
constraints) the two approaches are asymptotically equivalent. This
equivalence will be revisited in the present context.

\subsection{Score-based Indirect Inference With(out) Constraints}\label{score_ii}

 Given $H$ simulated paths $\{\tilde{y}_{t}^{(h)}(\theta )\}_{t=1}^{T},h=1,...,H$, a simulated version of the auxiliary criterion,
denoted by $Q_{TH}(\theta ,\beta )$, can then be constructed for use in I-I.
To fix ideas, say we have in mind auxiliary parameters $\beta $ defined as
M-estimators that maximize the criterion 
\[
Q_{T}(\beta )=\frac{1}{T}\sum_{t=1+l}^{T}q(y_{t},y_{t-1},..,y_{t-l};\beta ). 
\]%
The simulated auxiliary criterion $Q_{TH}(\theta ,\beta )$ is then
constructed by averaging over the $H$ paths\footnote{Note that our use of $Q_{TH}(\cdot)$ is a slight abuse of notation since, in the case of a dynamic model, the probability distribution of $Q_{TH}(\cdot)$ depends separately on $T$ and $H$ {and not only on the} product $TH$. This abuse of notation is immaterial for first-order asymptotics.}
\begin{equation}
{Q}_{TH}(\theta ,\beta )=\frac{1}{H}\sum_{h=1}^{H}\frac{1}{T}\sum_{t=l+1}^{T}q(\tilde{y}%
_{t}^{(h)}(\theta ),\tilde{y}_{t-1}^{(h)}(\theta ),..,y_{t-l}^{(h)}(\theta
);\beta ).  \label{THdef}
\end{equation}

Given ${Q}_{TH}(\theta ,\beta )$, under sufficient smoothness conditions, the
gradient (w.r.t. $\beta $) for the {simulated version of the quadratic criterion function $M_T(\beta)$} is given
by 
\[
\frac{\partial Q_{TH}(\theta ,{\beta }_{T}^{0})}{\partial \beta }+\frac{%
	\partial ^{2}Q_{TH}(\theta ,{\beta }_{T}^{0})}{\partial \beta \partial \beta
	^{\prime }}\left( {\beta }-{\beta }_{T}^{0}\right) . 
\]
Replacing the infeasible $\beta _{T}^{0}$ by $\hat{\beta}_{T}^{r}$, and
evaluating this gradient at $\beta =\widehat{\beta }_{T}$, we can then use
the resulting estimating equations
\begin{eqnarray}
\bar{m}_{TH}[\theta ;\widehat{\beta }_{T}] &=&\frac{\partial Q_{TH}(\theta ,%
	\hat{\beta}_{T}^{r})}{\partial \beta }+\frac{\partial ^{2}Q_{TH}(\theta ,%
	\hat{\beta}_{T}^{r})}{\partial \beta \partial \beta ^{\prime }}\left( 
\widehat{\beta }_{T}-\hat{\beta}_{T}^{r}\right)  \label{scoreus} \\
&=&\frac{\partial Q_{TH}(\theta ,\hat{\beta}_{T}^{r})}{\partial \beta }-%
\frac{\partial ^{2}Q_{TH}(\theta ,\hat{\beta}_{T}^{r})}{\partial \beta
	\partial \beta ^{\prime }}\left[ \frac{\partial ^{2}Q_{T}(\hat{\beta}%
	_{T}^{r})}{\partial \beta \partial \beta ^{\prime }}\right] ^{-1}\frac{%
	\partial Q_{T}(\hat{\beta}_{T}^{r})}{\partial \beta }  \nonumber
\end{eqnarray}%
to carry out a score-based I-I approach. In the absence of constraints for $%
\theta $, this approach yields the following I-I estimator 
\begin{equation}
\widehat{\theta }_{T,H}^{s}(W)=\arg \min_{\theta \in \Theta }\bar{m}%
_{TH}[\theta ;\widehat{\beta }_{T}]^{\prime }\cdot {W}\cdot \bar{m}%
_{TH}[\theta ;\widehat{\beta }_{T}],  \label{unr}
\end{equation}%
where $W$ is a positive-definite $(d_{\beta }\times d_{\beta })$ weighting
matrix.

In contrast to the I-I estimator in \eqref{unr},  the key idea of the CFS I-I strategy is to
incorporate the KT multipliers by considering the modified estimating equations\footnote{Note that CFS actually define this estimator only for ``$%
	H=\infty .$''}
\[
m_{TH}^{CFS}[\theta ;\hat{\lambda}_{T}]=\frac{\partial Q_{TH}(\theta ,\hat{%
		\beta}_{T}^{r})}{\partial \beta }+\frac{\partial g^{\prime }(\hat{\beta}%
	_{T}^{r})}{\partial \beta }\cdot \hat{\lambda}_{T} .
\]
However, by plugging in the vector $\hat{\lambda}_{T}$ of
KT multipliers, as given by the first-order conditions, \eqref{lagrange}, we obtain
\begin{equation}
m_{TH}^{CFS}[\theta ;\hat{\lambda}_{T}]=\frac{\partial Q_{TH}(\theta ,\hat{%
		\beta}_{T}^{r})}{\partial \beta }-\frac{\partial Q_{T}(\hat{\beta}_{T}^{r})}{%
	\partial \beta }  \label{CFS}.
\end{equation}
Then, for any
positive-definite $(d_{\beta }\times d_{\beta })$ weighting matrix $W$, CFS
 compute their so-called ``restricted" score-based I-I estimator as
\begin{equation}
\widehat{\theta }_{T,H}^{CFS}(W)=\arg \min_{\theta \in \Theta
}m_{TH}^{CFS}[\theta ;\hat{\lambda}_{T}]^{\prime}\cdot {W}\cdot
m_{TH}^{CFS}[\theta ;\hat{\lambda}_{T}]  \label{res}.
\end{equation}
CFS refer to their I-I estimator $\widehat{\theta }_{T,H}^{CFS}(W)$ as a restricted estimator, while we dub our I-I estimator $\widehat{\theta }_{T,H}^{s}(W)$ an ``unrestricted" estimator since we employ the unrestricted estimator $\widehat{\beta }_{T}$.\footnote{{We note here that both $\bar{m}_{TH}[\theta ;\widehat{\beta }_{T}]$ and $m_{TH}^{CFS}[\theta ;\hat{\lambda}_{T}]$ depend on the constrained estimator $\hat{\beta}^r_{T}$. While it is an abuse of notation to subsumed this dependence, we believe this avoids notational clutter and allows us to easily differentiate between the two estimating equations.} }

The first key result of this section is to demonstrate that the restricted terminology employed by CFS is potentially misleading: both the restricted equations used by CFS and our unrestricted equations are asymptotically equivalent to the estimating equations that would be used in an \textbf{\textit{unconstrained}} (but infeasible) GT-type score-based I-I approach. That is, the estimating equations used in both the CFS approach and our approach are equivalent to the following unrestricted estimating equations that match simulated data at the unconstrained, but infeasible, estimator $\ddot{\beta}_T$: 
\begin{equation}\label{GT}
\frac{\partial Q_{TH}(\theta,\ddot{\beta}_T)}{\partial\beta},\text{ where }\ddot{\beta}_{T}=\beta^0_T+J_{T}^{-1}\partial Q_T(\beta^0_T)/\partial\beta,
\end{equation}

Demonstrating equivalence between the three sets of estimating equations, \eqref{scoreus}, \eqref{CFS}, and \eqref{GT}, requires the following assumption.

\medskip

\noindent\textbf{Assumption A3:} Under $\{\theta_T\}\in\Gamma(\theta^0,\beta^0)$, the following are satisfied for any fixed $H\geq1$. 
\medskip 

\noindent(i) For all $\theta \in \Theta$, $Q_{TH}(\theta ,\beta )$ has
continuous partial derivatives (in $\beta$) of order two on $\text{Int}(\mathbf{B}%
_{{}}^{r})$ with probability one.
\medskip

\noindent(ii) For $\delta ^{0}$ and $\mathcal{I}^{0}$ defined in Assumption \textbf{A1}: 
\[
\sqrt{T}{\partial Q_{TH}(\theta }_{T}{,\beta _{T}^{0})}/{\partial \beta ^{{}}}%
\rightarrow _{d}\aleph \left( \delta ^{0},\mathcal{I}^{0}/H\right) . 
\]

\noindent(iii) There exists a continuous matrix function $\theta \mapsto \mathcal{J}(\theta
,\beta ^{0})$ such that, for all $\theta \in \Theta$,
$\mathcal{J}(\theta ,\beta ^{0})$ is positive-definite and for
any $\gamma >0:$ 
\[
\sup_{\theta \in \Theta }\sup_{\beta \in \mathbf{B}^{r}:\Vert \beta -\beta
	_{T}^{0}\Vert \leq \frac{\gamma }{\sqrt{T}}}\left\Vert \frac{\partial
	^{2}Q_{TH}(\theta ,\beta )}{\partial \beta \partial \beta ^{\prime }}+
\mathcal{J}(\theta ,\beta ^{0})\right\Vert =o_{P}(1). 
\]

While the contents of Assumption \textbf{A3} are relatively straightforward, we would like to point out that Assumption \textbf{A3}(ii) requires that, if we were to simulate
under the true value of the structural parameters,\ then asymptotically
the behavior of $\sqrt{T}{\partial Q_{TH}(\theta }_{T}{,\beta _{T}^{0})}/{%
	\partial \beta }$ and $\sqrt{T}\partial Q_{T}(\beta _{T}^{0})%
/{\partial \beta }$\ must agree. Implicitly, such an assumption
requires that the models we are simulating data from be
correctly specified, at least asymptotically.

Under Assumptions \textbf{A0}-\textbf{A3}, we have the following result.

\medskip 

\noindent\textbf{Proposition 2:} For $\{\theta_T\}\in\Gamma(\theta^0,\beta^0)$, under Assumptions \textbf{A0-A3}, for any given $H\geq 1$,
\begin{eqnarray*}
	\sup_{\theta \in \Theta }\left\Vert \bar{m}_{TH}[\theta ;\widehat{\beta }%
	_{T}]-{\partial Q_{TH}(\theta ,\ddot{\beta}_{T})}/{\partial \beta }%
	\right\Vert =o_{P}( 1/\sqrt{T})=\sup_{\theta \in \Theta }\left\Vert m_{TH}^{CFS}[\theta ;\hat{\lambda}_{T}]-%
	{\partial Q_{TH}(\theta ,\ddot{\beta}_{T})}/{\partial \beta }\right\Vert,
\end{eqnarray*}and it follows that $\sup_{\theta \in \Theta }\left\Vert m_{TH}^{CFS}[\theta ;\hat{\lambda}_{T}]-\bar{m}_{TH}[\theta ;\widehat{\beta }%
_{T}]\right\Vert
=o_{P}( 1/\sqrt{T}).$\hfill$\square$
\medskip 

From \textbf{Proposition 2}, we conclude that both the restricted estimating equations considered by CFS, and the unrestricted estimating equations proposed herein, are equivalent (uniformly over $\Theta$, at first-order) to unrestricted GT-type estimating equations. Therefore, the central message of \textbf{Proposition 2} is that constraints on the auxiliary model should have no effect on the choice of moments to match within score-based I-I estimation. 

To understand the significance of this asymptotic equivalence, recall that if the auxiliary parameters are near the boundary, the traditional approach based on minimizing, over $\Theta$, the estimating equations $\partial Q_{TH}(\theta,\hat{\beta}^r_T)/\partial\beta$ may not deliver an asymptotically Gaussian estimator of $\theta^0$. In contrast, an I-I approach based on the unconstrained (but infeasible) GT-type estimating equations in \eqref{GT} would deliver an asymptotically Gaussian estimator of $\theta^0$: from a Taylor series expansion, under Assumptions \textbf{A0}-\textbf{A3}, we can conclude that
\begin{flalign} \sqrt { T } \frac { \partial Q _ { T H }( \theta_{T} , \ddot { \beta } _ { T }) } { \partial \beta } & = \sqrt { T } \frac { \partial Q _ { T H } \left( \theta_{T} , \beta _ { T } ^ { 0 } \right) } { \partial \beta } + \frac { \partial ^ { 2 } Q _ { T H } \left( \theta_{T} , \beta^0_ { T } \right) } { \partial \beta \partial \beta ^ { \prime } } \sqrt { T } \left(\ddot { \beta } _ { T } - \beta _ { T } ^ { 0 } \right) +o_{P}(1)\nonumber\\ & = \sqrt { T } \frac { \partial Q _ { T H } \left( \theta_{T} , \beta _ { T } ^ { 0 } \right) } { \partial \beta } + \frac { \partial ^ { 2 } Q _ { T H } \left( \theta_{T} , \beta^0_ { T } \right) } { \partial \beta \partial \beta ^ { \prime } } \left\{ - \frac { \partial ^ { 2 } Q _ { T } \left( \beta _ { T } ^ { 0 } \right) } { \partial \beta \partial \beta ^ { \prime } } \right\} ^ { - 1 } \sqrt { T } \frac { \partial Q _ { T } \left( \beta _ { T } ^ { 0 } \right) } { \partial \beta }+o_{P}(1) \nonumber\\ & = \sqrt { T } \frac { \partial Q _ { T H } \left( \theta_{T} , \beta _ { T } ^ { 0 } \right) } { \partial \beta } - \sqrt { T } \frac { \partial Q _ { T } \left( \beta _ { T } ^ { 0 } \right) } { \partial \beta } + o _ { P } ( 1 ). \label{final}\end{flalign}
By Assumptions \textbf{A1}(iii) and \textbf{A3}(ii), {and the independence of the observed and simulated data,} the right-hand side term in equation \eqref{final} is asymptotically Gaussian with zero mean. Therefore, a direct consequence of \textbf{Proposition 2} and equation \eqref{final} is that our unrestricted estimating equations $\bar{m}_{TH}[\theta;\widehat{\beta}_T]$, when evaluated at $\theta_T$, are also asymptotically Gaussian even when the auxiliary parameters are near the boundary.

Given the asymptotic equivalence derived in \textbf{Proposition 2}, we would expect that our I-I estimator $\widehat{%
	\theta }_{T,H}^{s}(W)$ will be asymptotically equivalent to the CFS I-I
estimator $\widehat{\theta }_{T,H}^{CFS}(W)$. To detail such an equivalence result, we must maintain the following standard assumption for consistency of extremum estimators. 

\medskip 

\noindent\textbf{Assumption A4: }Under $\{\theta _{T}\}\in \Gamma (\theta ^{0},\beta
^{0})$, for any $H\ge1$:
\medskip

\noindent(i) For any $\beta \in $ Int$(\mathbf{B}^{r})$, the function $\theta \mapsto 
{\partial Q_{TH}(\theta ,\beta )}/{\partial \beta }$\ is continuous on $%
\Theta .$
\medskip

\noindent(ii) There exists a vector function $L(\theta ,\beta ^{0})$ such that, for
any $\gamma >0$, 
\[
\sup_{\theta \in \Theta }\sup_{\Vert \beta -\beta _{T}^{0}\Vert \leq \frac{%
		\gamma }{\sqrt{T}}}\left\Vert \frac{\partial Q_{TH}(\theta ,\beta )}{\partial
	\beta }-L(\theta ,\beta ^{0})\right\Vert =O_{P}\left( 1/\sqrt{T}\right) . 
\]

\noindent(iii) 
$
L\left( \theta ,\beta ^{0}\right) =0\Longleftrightarrow \theta =\theta ^{0}. 
$

\medskip

Assumption \textbf{A4} is a standard identification assumption, and is
implicitly maintained in CFS. However, our explicit treatment of parameters
near the boundary forces us to be more cautious. To see that, let us discuss the content
of the identification Assumption \textbf{A4} in the context of the
stochastic volatility model example in Section 2. For sake of expositional simplicity,
let us consider an auxiliary model based on conditional normality, with $%
\beta =(\psi ,\varphi ,\pi )^{\prime }.$ CFS rightly recall that $\pi $
becomes asymptotically underidentified when $\varphi =0.$ CFS circumvent
this issue by assuming $\varphi \geq 0.025.$ In contrast, we propose in this
paper an explicit treatment of parameters near the boundary, which may allow
the asymptotic true value $\beta ^{0}=(\psi ^{0},\varphi ^{0},\pi
^{0})^{\prime }$ \ to be such that $\varphi ^{0}=0$. The reader can easily
check that this specific value does not prevent ${\partial Q_{TH}(\theta
	,\beta ^{0})}/{\partial \beta }$ from having a well-defined probability limit 
$L(\theta ,\beta ^{0}).$ Consider a trial true value $\theta ^{0}=(\alpha
^{0},\delta ^{0},\sigma _{v}^{0})^{\prime }$ with $\delta ^{0}=0$. Then, $%
y_{t}$ is homoskedastic and 
\[
L(\theta ^{0},\beta ^{0})=0,
\]%
with $\beta ^{0}=(\psi ^{0},\varphi ^{0},\pi ^{0})^{\prime }$ \ and with 
\begin{eqnarray*}
	\psi ^{0} =\text{Var}(y_{t})=\alpha ^{0}, \text{ and }
	\varphi ^{0} =\pi ^{0}=0.
\end{eqnarray*}%
But, if $\theta =(\alpha ,\delta ,\sigma _{v})^{\prime }$ with $\delta \neq 0
$, then, $y_{t}$ is conditionally heteroskedastic and obviously:%
\[
L(\theta ,\beta ^{0})\neq 0.
\]%
From this toy example, we conclude that Assumptions \textbf{A4} is sensible.%
\footnote{%
	The reader may wonder how to identify $\sigma _{v}$ in the homoskedastic
	case. This actually requires matching the kurtosis since in the general case
	the unconditional kurtosis is 
	\[
	\frac{Var(h_{t})}{[E(h_{t})]^{2}}=\exp \left( \frac{\sigma _{v}^{2}}{%
		1-\delta ^{2}}\right) -1.
	\]%
	This kurtosis matching is implicitly performed when using a Student-t
	conditional distribution as an auxiliary model.}

Assumptions \textbf{A0}-\textbf{A4} allow us to prove
the following result.

\medskip

\noindent\textbf{Proposition 3:} For $\{\theta_T\}\in\Gamma(\theta^0,\beta^0)$, under Assumptions \textbf{A0-A4}, for any given $H\geq 1$ and any positive-definite
matrix $W$, $\plim_{T\rightarrow\infty}\widehat{\theta }_{T,H}^{s}(W)=\plim_{T\rightarrow\infty}\widehat{\theta }_{T,H}^{CFS}(W)=\theta^0$, and
$\left\| \widehat{\theta }_{T,H}^{s}(W)-\widehat{\theta }%
_{T,H}^{CFS}(W)\right\| =o_{P}(1/\sqrt{T})$.
\hfill$\square$

\medskip 

 Since our unrestricted I-I estimator $\widehat{\theta }_{T,H}^{s}(W)$ is
 asymptotically equivalent to the restricted I-I estimator $\widehat{\theta }%
 _{T,H}^{CFS}(W)$, we will set the focus on the former. By doing so, we
 confirm the discussion {given earlier and in Section \ref{sec2}:} when it comes to
 the choice of the moments to match, we do not really care {about} constrained
 estimation of the auxiliary model.

To prove asymptotic normality of the resulting estimator, a local
identification assumption is required to complete the global Assumption \textbf{A4}(iii). \medskip
 
 \noindent \textbf{Assumption A5}: {The vector function $\theta \mapsto
 L(\theta ,\beta ^{0})$ is continuously differentiable on $\text{Int}(\Theta )
 $\ and, for all $T\geq1$,  $
 \text{rank}\left( {\partial L(\theta _{T},\beta ^{0})}/{\partial \theta
 ^{\prime }}\right) =d_{\theta }= \text{rank}\left( {\partial L(\theta ^{0},\beta ^{0})}/{\partial \theta
 ^{\prime }}\right).
$}
 
 \medskip 
 
 \noindent\textbf{Theorem 2}: For $\{\theta_T\}\in\Gamma(\theta^0,\beta^0)$, under Assumptions \textbf{A0-A5}, for any given $H\geq 1$ and any positive-definite
 matrix $W$
 \[
 \sqrt{T}\left( \widehat{\theta }_{T,H}^{s}(W)-\theta_{T}\right) \rightarrow
 _{d}\aleph \left( 0,\left( 1+\frac{1}{H}\right) \Omega _{W}\right) ,
 \]
 where, recalling $\mathcal{I}^0=\lim_{T\rightarrow\infty}\text{Var}\left[\sqrt{T}{\partial Q_{T}(\beta^0_T)}/{\partial\beta}\right]$, 
 \begin{eqnarray*}
\Omega _{W}=A_{W}^{-1}B_{W}A_{W}^{-1},\;\;  A_{W} =\frac{\partial L(\theta ^{0},\beta ^{0})^{\prime }}{\partial \theta 
 }W\frac{\partial L(\theta ^{0},\beta ^{0})}{\partial \theta ^{\prime }}, \;\;
 B_{W} =\frac{\partial L(\theta ^{0},\beta ^{0})^{\prime }}{\partial \theta 
 }W\mathcal{I}^{0}W\frac{\partial L(\theta ^{0},\beta ^{0})}{\partial \theta
 ^{\prime }}.
 \end{eqnarray*}
The optimal weighting matrix $W$ is given by%
 \[
 W^{\ast }=[\mathcal{I}^{0}]^{-1},
  \]
 and leads to an optimal I-I estimator with asymptotic variance\footnote{%
 The reader may notice that the formula for $\Omega ^{*}$ given above differs
 from that given in {Proposition 4} of CFS, denoted as $\mathcal{C}_{0}^{r}$
 in their equation (8). However, it is simple to verify that the two coincide
 when the constraints $g(\beta)\geq0$ are satisfied at $\beta^0_T$ since the KT multipliers will be zero in the limit.}%
 \[
 \left( 1+\frac{1}{H}\right) \Omega ^{\ast }=\left( 1+\frac{1}{H}\right)
 \left( \frac{\partial L(\theta ^{0},\beta ^{0})^{\prime }}{\partial \theta }[%
 \mathcal{I}^{0}]^{-1}\frac{\partial L(\theta ^{0},\beta ^{0})}{\partial
 \theta ^{\prime }}\right) ^{-1}.
 \]
 \hfill$\square$

\medskip 

{Theorem 2 demonstrates that even though, due to the boundary bias in Assumption \textbf{A2}, the term $\sqrt{T}(\theta_T-\theta^0)$ may not converge to zero, the term} $\sqrt{T}( \widehat{\theta }_{T,H}^{s}(W)-\theta _{T}) $ still converges to a zero-mean Gaussian random variable. This result comes about from the structure of $\bar{m}_{TH}[\theta;\widehat{\beta}_{T}]$ and the results of \textbf{Proposition 2}, which imply that
$$
\sqrt{T}\bar{m}_{TH}[\theta_T;\widehat{\beta}_{T}]=\sqrt{T}\frac { \partial Q _ { T H } \left( \theta_T ,  { \beta } _ { T } ^ { 0 } \right) } { \partial \beta } -\sqrt{T}\frac { \partial Q _ { T } \left( { \beta } _ { T } ^ { 0 } \right) } { \partial \beta ^ { } }+o_{P}(1).
$$Then, Assumptions \textbf{A1}(iii) and \textbf{A3}(ii), together with the independence between the observed and simulated data, ensure that $\{\sqrt{T}\partial Q_{TH}(\theta_T,\beta^0_T)/\partial\beta-\sqrt{T}\partial Q_{T}(\beta^0_T)/\partial\beta\}$ is asymptotically Gaussian with zero mean. That is, even though $ \sqrt{T}{ \partial Q _ { T } \left( { \beta } _ { T } ^ { 0 } \right) }/  \partial \beta ^ {  }$ has a non-zero asymptotic mean, because I-I seeks to simulate data so that the observed and simulated scores agree, in the sense that their normed difference is small, this non-zero asymptotic mean is ``knocked out'' and does not contaminate the asymptotic distribution of $\sqrt{T}(\hat{\theta}_{T,H}^{s}(W)-\theta_T)$.

It is also important to note that the above formulas are identical to those given in
GMR, confirming that we actually perform I-I \textit{\textbf{without}} constraints. To see this, let $Q_{T}(\theta
,\beta )$ denote the simulated auxiliary criterion function calculated using the single simulated path ($H=1$) $\{\tilde{y}_{t}(\theta )\}_{t=1}^{T}$, i.e., with reference to equation \eqref{THdef}
$$
Q_{T}(\theta,\beta ):=\frac{1}{T}\sum_{t=l+1}^{T}q(\tilde{y}%
_{t}^{(h)}(\theta ),\tilde{y}_{t-1}^{(h)}(\theta ),..,y_{t-l}^{(h)}(\theta
);\beta ),
$$and consider the constrained estimator
\[
\tilde{\beta}_{T}^{r}(\theta )=\arg \max_{\beta \in \mathbf{B}^{r}}Q_{T}(\theta
,\beta ).
\]
For sake of interpretation, let us consider the simplest case without
boundary problems. Then, the constrained estimator $\tilde{\beta}%
_{T}^{r}(\theta )$\ converges towards a (non-drifting) pseudo-true value $%
b(\theta )$ that is in the interior of the parameter set. Then, while KT
multipliers converge to zero, we have
\[
\Plim{T\rightarrow\infty }\frac{\partial Q_{T}( \theta, \tilde{\beta}%
_{T}^{r}(\theta ))}{\partial \beta } =L\left( \theta ,b(\theta )\right) =0,\forall
\theta \in \Theta .
\]
{In particular, by differentiating the above and assuming, following Assumption \textbf{A1}, $$\frac{\partial L(\theta^0,\beta^0)}{\partial\beta'}=\plim_{T\rightarrow\infty} \frac{\partial^2 Q_T(\beta^0)}{\partial\beta\partial\beta'}=-\mathcal{J}^0
$$we obtain
\[
\frac{\partial L(\theta ^{0},\beta ^{0})}{\partial \theta ^{\prime }}+\frac{%
\partial L(\theta ^{0},\beta ^{0})}{\partial \beta ^{\prime }}\frac{\partial
b(\theta ^{0})}{\partial \theta ^{\prime }}=0
\]
and so}
\[
\frac{\partial L(\theta ^{0},\beta ^{0})}{\partial \theta ^{\prime }}=\mathcal{J}^{0}\frac{\partial b(\theta ^{0})}{\partial \theta ^{\prime }}.
\]
{This relationship between ${\partial L(\theta ^{0},\beta ^{0})}/{\partial \theta ^{\prime }}$ and $\partial b(\theta)/\partial\theta'$ allows us to rewrite the asymptotic variance of $\sqrt{T}(\widehat{\theta}^{s}_{T,H}(W)-\theta_T)$ as }
\begin{flalign*}
\Omega_{W} &= A_{W}^{-1}B_{W} A_{W}^{-1},\;\;A_{W}=\frac{\partial b(\theta^0)'}{\partial\theta}\mathcal{J}^{0'}W\mathcal{J}^{0}\frac{\partial b(\theta^0)}{\partial\theta'},\;\;
B_{W} =\frac{\partial b(\theta^0)'}{\partial\theta}\mathcal{J}^{0'}W\mathcal{I} W\mathcal{J}^{0}\frac{\partial b(\theta^0)}{\partial\theta'}.
\end{flalign*}
Therefore,
\[
\Omega ^{\ast }=\left\{ \frac{\partial b^{\prime }(\theta ^{0})}{\partial
\theta }\mathcal{J}^{0'}[\mathcal{I}^{0}]^{-1}\mathcal{J}^{0}\frac{\partial b(\theta ^{0})}{\partial
\theta ^{\prime }}\right\} ^{-1},
\]
and we recognize the familiar formula given by GMR (see their Proposition 4) for the
asymptotic variance of the optimal I-I estimator.

\subsection{Wald-based Indirect Inference With(Out) Constraints}

 The aforementioned tight connection with the results of GMR
suggest that it should be possible to perform I-I \textit{\textbf{without}}
constraints in an alternative, albeit asymptotically equivalent, {manner using the} Wald approach {and our well-behaved unconstrained estimator $\widehat{\beta}_{T}$}. {The philosophy} of the Wald
approach to I-I would {then} amount to compute an unconstrained estimator $%
\tilde{\beta}_{TH}(\theta )$ on simulated data (for any given value $\theta $%
\ of the structural parameters) and then to minimize{, in some norm,} $\widehat{\beta}_{T}-\tilde{\beta}_{TH}(\theta )$. We show in this section that this approach may work, but
{requires} care in the definition of $\tilde{\beta}_{TH}(\theta )$ .

\subsubsection{A First Solution: the CFS Strategy}

The {Wald-based I-I strategy of CFS}, {which uses constrained auxiliary parameter estimates,} can be reinterpreted as a
minimum distance I-I approach {based on a vector of unconstrained auxiliary parameter estimates}.
To see this, first define 
\[
\tilde{\beta}_{TH}^{r}(\theta )=\arg \max_{\beta \in \mathbf{B}^{r}}{Q}%
_{TH}(\theta ,\beta ),
\]
and let $\tilde{\lambda}_{TH}(\theta )$ be the vector of KT multipliers
delivered by this constrained optimization. The Wald-based estimator of CFS is then given by
\begin{equation}
\check{\theta}_{T,H}^{CFS}(W)=\arg \min_{\theta \in \Theta }\left[ 
\begin{array}{c}
\hat{\beta}_{T}^{r}-\tilde{\beta}_{TH}^{r}(\theta ) \\ 
\hat{\lambda}_{T}-\tilde{\lambda}_{TH}(\theta )%
\end{array}%
\right] ^{\prime }K_{0}^{r\prime }\cdot {W}^{\boxplus }\cdot K_{0}^{r}\left[ 
\begin{array}{c}
\hat{\beta}_{T}^{r}-\tilde{\beta}_{TH}^{r}(\theta ) \\ 
\hat{\lambda}_{T}-\tilde{\lambda}_{TH}(\theta )%
\end{array}%
\right],   \label{infeaW}
\end{equation}
where, {since we are under the assumption} that the constraints are fulfilled, 
\[
{W}^{\boxplus }=\left[ 
\begin{array}{cc}
W & \mathbf{O} \\ 
\mathbf{O} & \mathbf{O}%
\end{array}%
\right],\;\;K^r_0=\begin{bmatrix}K_{0,1}^r\\\mathbf{O}\end{bmatrix},\;\; 
K_{0,1}^{r}=\left[ 
\begin{array}{ccc}
	-\mathcal{J}^{0} & \vdots &\frac{\partial g^{\prime }(\beta _{T}^{0})}{\partial \beta 
	}%
\end{array}%
\right].
\] 

 {Recall that we
have simplified} the exposition by considering {only auxiliary
parameter estimates $\tilde{\beta}_{TH}^{r}(\theta )$ defined as above}. {{}Alternatively, we could consider $H$ auxiliary parameters based on a single simulated paths of length $T$: for $h=1,\dots,H$}
\begin{eqnarray*}
\tilde{\beta}_{T}^{r(h)}(\theta ) &=&\arg \max_{\beta \in \mathbf{B}^{r}}{{Q}}_{T}^{}(\theta ,\beta ),
\end{eqnarray*}and then compute\footnote{Extending the results of GMR, we can conclude that $\tilde{%
\beta}_{TH}^{r}(\theta )$ and $\bar{\beta}_{T,H}^{r}(\theta )$ are
asymptotically equivalent and would lead to asymptotically equivalent I-I estimators of $\theta$. {However,}
the results of Gourieroux, Renault and Touzi (2000) {suggest} that an I-I
estimator based on $\bar{\beta}_{T,H}^{r}(\theta )$ will have better
finite sample properties, at the cost of performing $H$ optimizations in the
auxiliary model instead of only just one. This discussion is beyond the
scope of this paper.
}
 $$
\bar{\beta}_{T,H}^{r}(\theta ) =\frac{1}{H}\sum_{h=1}^{H}\tilde{\beta}%
_{T}^{(h)}(\theta ).$$

{Note that the estimator $\bar{\beta}_{T,H}^{r}(\theta )$ fulfills the constraints if $\mathbf{B}^{r}$\ is a
convex set. A sufficient condition for that is to assume that the set $\mathbf{B}_{}$ %
is convex and the functions $g_{j}(\cdot),j=1,..,q$ defining the constraints are concave. However, the satisfaction of this condition is immaterial for the validity of an I-I estimator based on $\bar{\beta}_{T,H}^{r}(\theta )$. }

{It must be acknowledged that a}
more complicated definition for $K_{0,1}^{r}$ {is given in CFS. However,} this
complication {is immaterial in our setting as we work} {under the assumption} that the constraints are fulfilled, and thus the population
(resp., estimated) vector of KT multipliers is {zero} (resp., $O_{P}( 1/\sqrt{T}) $). As a matter of fact, the above estimator becomes feasible only when 
$K_{0,1}^{r}$ is replaced by a consistent estimator like
\[
\hat{K}_{0,1,T}^{r}=\left[ 
\begin{array}{ccc}
\frac{\partial ^{2}Q_{T}(\hat{\beta}_{T}^{r})}{\partial \beta
\partial \beta ^{\prime }} &\vdots& \frac{\partial g^{\prime }(\hat{\beta}_{T}^{r})%
}{\partial \beta }%
\end{array}%
\right] =\left[ 
\begin{array}{ccc}
-\hat{J}_{T} &\vdots& \frac{\partial g^{\prime }(\hat{\beta}_{T}^{r})}{\partial
\beta }%
\end{array}%
\right] .
\]

Hence, for the sake of feasibility, we should rather consider
\begin{equation}
\check{\theta}_{T,H}^{CFS}(W)=\arg \min_{\theta \in \Theta }\left[ 
\begin{array}{c}
\hat{\beta}_{T}^{r}-\tilde{\beta}_{TH}^{r}(\theta ) \\ 
\hat{\lambda}_{T}-\tilde{\lambda}_{TH}(\theta )%
\end{array}%
\right] ^{\prime }\hat{K}_{0,1,T}^{r\prime }\cdot {W}\cdot \hat{K}%
_{0,1,T}^{r}\left[ 
\begin{array}{c}
\hat{\beta}_{T}^{r}-\tilde{\beta}_{TH}^{r}(\theta ) \\ 
\hat{\lambda}_{T}-\tilde{\lambda}_{TH}(\theta )%
\end{array}%
\right] .  \label{feaW}
\end{equation}
 Since the two estimators (\ref{infeaW}) and (\ref{feaW}) are
obviously asymptotically equivalent, we simplify the exposition by denoting
them identically, even though only (\ref{feaW}) is feasible. 

{Just as with the score-based approach to I-I}, {we can now} interpret the
Wald-based I-I estimator {of CFS} as I-I \textit{\textbf{with{out}}} constraints. {To do so, note that}
\begin{eqnarray*}
\hat{K}_{0,1,T}^{r}\left[ 
\begin{array}{c}
\hat{\beta}_{T}^{r}-\beta _{T}^{0} \\ 
\hat{\lambda}_{T}%
\end{array}%
\right]  &=&\frac{\partial ^{2}Q_{T}(\hat{\beta}_{T}^{r})}{\partial \beta
\partial \beta ^{\prime }}\left( \hat{\beta}_{T}^{r}-\beta _{T}^{0}\right) +%
\frac{\partial g^{\prime }(\hat{\beta}_{T}^{r})}{\partial \beta }\hat{\lambda%
}_{T} \\
&=&-\hat{J}_{T}\left( \widehat{\beta}_{T}-\beta _{T}^{0}\right) +o_{P}\left( 1/\sqrt{T}%
\right) ,
\end{eqnarray*}
where the second equality follows from equation \eqref{13++}. Therefore, an asymptotically
equivalent version of the CFS Wald-based I-I estimator could be computed as
\[
\bar{\theta}_{T,H}^{CFS}(W)=\arg \min_{\theta \in \Theta }\left( \widehat{\beta}%
_{T}-\tilde{\beta}_{TH}^{CFS}(\theta )\right) ^{\prime }\hat{J}_{T}W\hat{J}_{T}\left( 
\widehat{\beta}_{T}-\tilde{\beta}_{TH}^{CFS}(\theta )\right) ,
\]
where we define $\tilde{\beta}_{TH}^{CFS}(\theta )$ as
\begin{equation}
\tilde{\beta}_{TH}^{CFS}(\theta )=\tilde{\beta}_{TH}^{r}(\theta )+\left[ 
\frac{\partial ^{2}Q_{T}(\hat{\beta}_{T}^{r})}{\partial \beta \partial \beta
^{\prime }}\right] ^{-1}\frac{\partial g^{\prime }(\hat{\beta}_{T}^{r})}{%
\partial \beta }\tilde{\lambda}_{TH}(\theta ).
\label{WII1}
\end{equation}
Note that the notation $\tilde{\beta}_{TH}^{CFS}(\theta )$ is
justified by analogy with the relationships
\begin{eqnarray}
\widehat{\beta }_{T} &=&\hat{\beta}_{T}^{r}-\left[ \frac{\partial ^{2}Q_{T}(%
\hat{\beta}_{T}^{r})}{\partial \beta \partial \beta ^{\prime }}\right] ^{-1}%
\frac{\partial Q_{T}(\hat{\beta}_{T}^{r})}{\partial \beta ^{ }}
\label{FUNC1} \\
&=&\hat{\beta}_{T}^{r}+\left[ \frac{\partial ^{2}Q_{T}(\hat{\beta}_{T}^{r})}{%
\partial \beta \partial \beta ^{\prime }}\right] ^{-1}\frac{\partial
g^{\prime }(\hat{\beta}_{T}^{r})}{\partial \beta }\hat{\lambda}_{T}.
\label{FUNC2}
\end{eqnarray}

This reinterpretation of the so-called ``restricted'' Wald approach to I-I, as
dubbed by CFS, is an unconstrained I-I approach based (through equation (\ref{feaW})) on our FUNC
estimator. Therefore, we have a similar message to the score-based
approach. This is confirmed by {Proposition 5} and {6} of CFS,
which yield the following insights.

\medskip 

\noindent{(i)} For any choice of the positive definite weighting matrix $W$ (or more
generally for {any sequence} of sample dependent positive-definite
weighting matrices $W_{T}$ with a positive-definite limit), the score-based
I-I estimator $\widehat{\theta}_{T,H}^{CFS}(W)$  and the Wald-based I-I
estimator $\check{\theta}_{T,H}^{CFS}(W)$\ are asymptotically equivalent.
\medskip 

\noindent{(ii)} For $T$ sufficiently large, the two estimators are numerically equal
in the case of an auxiliary model that just identifies the structural
parameters because $d_{\beta }=d_{\theta }$.
\medskip 

Point {(i)} above revisits the results of GMR (see
their Section 2.5 page S91), {demonstrating that,} for any choice of the weighting
matrix $W$, the score-based approach with weighting matrix $W$ is
asymptotically equivalent to the Wald-based approach with weighting matrix $%
\hat{J}_{T}W\hat{J}_{T}$ as in the definition of $\widehat{\theta}^{CFS}_{T}(W)$. In the {case of a just identified auxiliary model,}
the choice of the weighting matrix is immaterial and point
(ii) {calls to mind} {Proposition 4.1}. in Gourieroux and
Monfort (1996). Once more, this {similarity to the} results of GMR and Gourieroux
and Monfort (1996) {confirms} that we are actually
performing I-I \textit{\textbf{without}} constraints. {In addition,} since our unrestricted
score-based I-I estimator $\widehat{\theta}_{T,H}^{s}(W)$ is
asymptotically equivalent to the restricted score-based estimator $\widehat{\theta }%
_{T,H}^{CFS}(W)$ (see {Theorem 2}), it is also (by point (i) above) asymptotically
equivalent to the alternative aforementioned Wald-based {estimators of CFS:} $%
\check{\theta}_{T,H}^{CFS}(W)$ and $\bar{\theta}_{T,H}^{CFS}(W)$\ .

\subsubsection{A Second Solution: Back to the Score}

The previous subsection revisited} the Wald-based CFS
estimator by resorting to a definition of $\tilde{\beta}_{TH}(\theta )$
that mimics, on simulated data, the {alternative definition of the FUNC estimator given in equation \eqref{FUNC2}}. We can alternatively use a definition of $\tilde{\beta}%
_{TH}(\theta )$ that mimics equation (\ref{FUNC1}). To see this,
recall that our score-based approach was focused on minimizing, in some norm,
\begin{eqnarray*}
\bar{m}_{TH}[\theta ;\widehat{\beta }_{T}] &=&\frac{\partial Q_{TH}(\theta ,%
\hat{\beta}_{T}^{r})}{\partial \beta }+\frac{\partial ^{2}Q_{TH}(\theta ,%
\hat{\beta}_{T}^{r})}{\partial \beta \partial \beta ^{\prime }}\left( 
\widehat{\beta }_{T}-\hat{\beta}_{T}^{r}\right)  \\
&=&\left[ \frac{\partial ^{2}Q_{TH}(\theta ,\hat{\beta}_{T}^{r})}{\partial
\beta \partial \beta ^{\prime }}\right] \left\{ \widehat{\beta }_{T}-\hat{%
\beta}_{T}^{r}+\left[ \frac{\partial ^{2}Q_{TH}(\theta ,\hat{\beta}_{T}^{r})%
}{\partial \beta \partial \beta ^{\prime }}\right] ^{-1}\frac{\partial
Q_{TH}(\theta ,\hat{\beta}_{T}^{r})}{\partial \beta }\right\}  \\
&=&\left[ \frac{\partial ^{2}Q_{TH}(\theta ,\hat{\beta}_{T}^{r})}{\partial
\beta \partial \beta ^{\prime }}\right] \left\{ \widehat{\beta }_{T}-\tilde{%
\beta}_{T,H}^{c}(\theta )\right\} 
\end{eqnarray*}where
\begin{equation}
\tilde{\beta}_{TH}^{c}(\theta )=\hat{\beta}_{T}^{r}-\left[ \frac{\partial
^{2}Q_{TH}(\theta ,\hat{\beta}_{T}^{r})}{\partial \beta \partial \beta
^{\prime }}\right] ^{-1}\frac{\partial Q_{TH}(\theta ,\hat{\beta}_{T}^{r})}{%
\partial \beta ^{}}  \label{cheat}.
\end{equation}

Let us acknowledge, however, an important difference of philosophy between the
definitions of $\tilde{\beta}_{TH}^{CFS}(\theta )$ and $\tilde{\beta}%
_{TH}^{c}(\theta )$. In the former case, we make a Newton-Raphson
improvement of $\tilde{\beta}_{TH}^{r}(\theta )$, while in the latter case
we remain true to $\hat{\beta}_{T}^{r}$. In this respect, we obviously set
the focus on score matching and, as a consequence, a comparison with our
score-based approach is straightforward. More precisely, if we define another Wald-based I-I estimator, the solution of
\[
\widehat{\theta}_{T,H}^{c}(W)=\arg \min_{\theta \in \Theta }\left( \widehat{%
\beta }_{T}-\tilde{\beta}_{TH}^{c}(\theta )\right) ^{\prime }\hat{J}_{T}W%
\hat{J}_{T}\left( \widehat{\beta }_{T}-\tilde{\beta}_{TH}^{c}(\theta
)\right) ,
\]
we see that, from the formulas above, this minimization program can be equivalently written as
\begin{flalign}
\min_{\theta \in \Theta }\bar{m}_{TH}[\theta ;\widehat{\beta }_{T}]\left[ \frac{%
\partial ^{2}Q_{TH}(\theta ,\hat{\beta}_{T}^{r})}{\partial \beta \partial
\beta ^{\prime }}\right] ^{-1}\hat{J}_{T}W\hat{J}_{T}\left[ \frac{\partial
^{2}Q_{TH}(\theta ,\hat{\beta}_{T}^{r})}{\partial \beta \partial \beta
^{\prime }}\right] ^{-1}\bar{m}_{TH}[\theta ;\widehat{\beta }_{T}] , \label{inter} 
\end{flalign}
which is nothing but minimizing a certain norm of $\bar{m}_{TH}[\theta ;\widehat{\beta }_{T}]$ exactly as in equation \eqref{unr}.\footnote{ It might be argued that we are not exactly minimizing a
norm w.r.t. $\theta $ since the weighting matrix itself depends on $\theta $%
. However, it must be realized that this is immaterial, both for consistency
and asymptotic distribution, to replace the occurrence of $\theta $\ in the
weighting matrix\ by a first-step consistent estimator. This argument is
quite similar to the one of equivalence between continuously updated GMM
(Hansen, Heaton and Yaron, 1996) and efficient two-step GMM.} As observed by a referee, for the above Wald-based I-I estimator, the auxiliary parameters are never estimated on simulated data, only on the observed data, which can be convenient if the chosen auxiliary model is computationally challenging in some manner{, and signifies that we are basically back to a ``score-based'' I-I approach. }

The asymptotic distribution of $\widehat{\theta}_{T,H}^{c}(W)$\
obviously depends on the limit of the {weighting matrix sequence given by}
\[
\Plim{T\rightarrow\infty }\left[ \frac{\partial ^{2}Q_{TH}(\theta_{T} ,\hat{\beta}%
_{T}^{r})}{\partial \beta \partial \beta ^{\prime }}\right] ^{-1}\hat{J}_{T}W%
\hat{J}_{T}\left[ \frac{\partial ^{2}Q_{TH}(\theta_{T} ,\hat{\beta}_{T}^{r})}{%
\partial \beta \partial \beta ^{\prime }}\right] ^{-1}=W.
\]
We can then conclude that this Wald-based I-I estimator $\widehat{\theta%
}_{T,H}^{c}(W)$ is asymptotically equivalent to the score-based I-I
estimator $\widehat{\theta}_{T,H}^{s}(W)$ introduced in Subsection \ref{score_ii}.
 In other words, all I-I estimators discussed so far (for the same
weighting matrix $W$) are asymptotically equivalent, exactly as in GMR. Furthermore, as in Theorem 2 above, the optimal choice of $W$ is $W^{\ast}=[\mathcal{I}^0]^{-1}$.

{ Interestingly enough, our unconstrained view of I-I results in numerical equivalence between this
Wald-based I-I estimator and our score-based I-I estimator when the dimension of the auxiliary and structural parameters are equal.}

\medskip 

\noindent\textbf{Theorem 3:}
For $T$\ sufficiently large and in the case of
a just identified auxiliary model ($%
d_{\beta }=d_{\theta }$), the estimators $\widehat{\theta}_{T,H}^{s}(W)$\
and $\widehat{\theta}_{T,H}^{c}(W^{\ast })$ are numerically identical irrespective of the choice of weighting matrix (i.e., $W\neq W^{*}$).\hfill$\square$

\medskip 

 To conclude this subsection, it is worth comparing, in more detail, the two
definitions of $\tilde{\beta}_{TH}(\theta )$ that have delivered Wald-based
I-I estimators (by calibration against the FUNC estimator) that are
asymptotically equivalent to the score-based approach:%
\begin{flalign*}
\tilde{\beta}_{TH}^{CFS}(\theta )&=\tilde{\beta}_{TH}^{r}(\theta )+\left[ 
\frac{\partial ^{2}Q_{T}(\hat{\beta}_{T}^{r})}{\partial \beta \partial \beta
^{\prime }}\right] ^{-1}\frac{\partial g^{\prime }(\hat{\beta}_{T}^{r})}{%
\partial \beta }\tilde{\lambda}_{TH}(\theta ),\\
\tilde{\beta}_{TH}^{c}(\theta )&=\hat{\beta}_{T}^{r}-\left[ \frac{\partial
^{2}Q_{TH}(\theta ,\hat{\beta}_{T}^{r})}{\partial \beta \partial \beta
^{\prime }}\right] ^{-1}\frac{\partial Q_{TH}(\theta ,\hat{\beta}_{T}^{r})}{%
\partial \beta}.
\end{flalign*}
Since $\tilde{\beta}_{TH}^{CFS}(\theta )$ is based on constrained
estimation on the simulated path, through the computation of $\tilde{\beta}%
_{TH}^{r}(\theta )$\ and $\tilde{\lambda}_{TH}(\theta )$, one may {wish to revisit} $\tilde{\beta}_{TH}^{c}(\theta )$ by also using
constrained estimators on the simulated path, that is by instead computing
\begin{eqnarray*}
\tilde{\beta}_{TH}^{func}(\theta ) &=&\tilde{\beta}_{TH}^{r}(\theta )-%
\left[ \frac{\partial ^{2}Q_{TH}(\theta ,\tilde{\beta}_{TH}^{r}(\theta ))}{%
\partial \beta \partial \beta ^{\prime }}\right] ^{-1}\frac{\partial
Q_{TH}(\theta ,\tilde{\beta}_{TH}^{r}(\theta ))}{\partial \beta ^{ }}
,\\
&=&\tilde{\beta}_{TH}^{r}(\theta )+\left[ \frac{\partial ^{2}Q_{TH}(\theta ,%
\tilde{\beta}_{TH}^{r}(\theta ))}{\partial \beta \partial \beta ^{\prime }}%
\right] ^{-1}\frac{\partial g^{\prime }(\tilde{\beta}_{TH}^{r}(\theta ))}{%
\partial \beta }\tilde{\lambda}_{TH}(\theta ).
\end{eqnarray*}

$\tilde{\beta}_{TH}^{func}(\theta )$ is the FUNC
estimator computed on the simulated path, and it seems sensible to match it
against the FUNC estimator $\widehat{\beta }_{T}$ computed on the observed
data. However, this approach will not deliver a
consistent estimator of $\theta^{0}$ in general. To see this, note that $\tilde{\beta}%
_{TH}^{CFS}(\theta )$ and\ $\tilde{\beta}_{TH}^{func}(\theta )$ both
set the focus on the same linear combination of $\tilde{\beta}%
_{TH}^{r}(\theta )$\ and $\tilde{\lambda}_{TH}(\theta )$. However, while $%
\tilde{\beta}_{TH}^{CFS}(\theta )$\ is guaranteed to end up with a
consistent estimator for the coefficients of this linear combination, the
coefficients in $\tilde{\beta}_{TH}^{func}(\theta )$\ themselves depend on
the unknown $\theta $. {As a consequence, setting the focus on $[\widehat{\beta}_{T}-\tilde{\beta}_{TH}^{func}(\theta )]$ alone, or some norm thereof, can induce an additional perverse solution in the {limit; i.e., }the limiting estimating equations, because of their nonlinear dependence on $\theta$, can admit an additional solution $\bar{\theta}$ with $\bar{\theta}\neq\theta^0$. As such, an I-I strategy based on $\tilde{\beta}_{TH}^{func}(\theta )$ above may not identify $\theta^{0}$.\footnote{Frazier and Renault (2017) give additional examples of settings where such perverse roots can arise in nonlinear econometric models. }}

\section{Illustrative Examples}
\subsection{Stochastic Volatility}\label{sv_ex}
In this section, we apply our score-based I-I approach to estimate the parameters of the stochastic volatility (SV) model:
\begin{flalign}
y_t& = \sqrt{h_{t}}e_{t}\label{vol1},\\\ln(h_{t})&=\alpha +\delta \ln(h_{t-1})+\sigma_{v}v_{t},\label{vol2}
\end{flalign}where $0<\delta<1$, $\sigma_{v}>0$, $(e_{t},v_{t})^{\prime}\sim_{i.i.d.}N(0,\text{Id}_{2})$ and $\theta=(\alpha,\delta,\sigma_{v})^{\prime}$. We observe a series $\{y_{t}\}_{t=1}^{T}$ from the SV model in \eqref{vol1}-\eqref{vol2} and our goal is to conduct inference on $\theta$.

Following the discussion in Section \ref{sec2}, we consider the GARCH(1,1) auxiliary model
\begin{flalign}
y_t& = \sqrt{h_{t}}\epsilon_{t}\label{auxvol1}\\ h_{t}&=\psi +\varphi y_{t-1}^{2} +\pi h_{t-1}^{2}\nonumber
\end{flalign} where the errors $\epsilon_{t}$ in \eqref{auxvol1} are $\epsilon_{t}\sim_{iid} N(0,1)$. The auxiliary parameters are denoted by $\beta$, with $\beta=(\psi,\varphi,\pi)^{\prime}$. As mentioned in Section \ref{sec2}, to ensure the GARCH(1,1) auxiliary model is well-behaved CFS require the following inequality constraint
\begin{eqnarray*}
\varphi&\geq.025.
\end{eqnarray*}
Unlike the approach of CFS, by considering drifting sequences of auxiliary parameters, we allow the constrained estimator to fully reach the boundary of constrained space, in the limit. That is, instead, we assume the true auxiliary parameters satisfy the inequality
\begin{eqnarray*}
\varphi_{T}^0&\geq&o(1).
\end{eqnarray*}

\subsubsection{Monte Carlo Design}
To assess the performance of our proposed I-I estimation strategy we follow the Monte Carlo design of Jacquier, Polson and Rossi (1994) (JPR, hereafter), also used in CFS. In particular, we consider two sets of structural parameters: $\theta^{0,1} = (-.736,.90,.363)^{\prime}$ and $\theta^{0,2} = (-.147,.98,.0614)^{\prime}$. These particular values for $\theta^0$ are related to the unconditional coefficient of variation $\kappa$ for the unobserved level of volatility $h_{t}$, where $$\kappa^2 = \frac{\text{Var}(h_{t})}{\left(E[h_{t}]\right)^2}=\exp\left(\frac{\sigma^2_v}{1-\delta}\right) -1.$$In the first design,  we have $\kappa^2=1$, which roughly represents lower-frequency returns (say, weekly or monthly returns); for the second design, we set $\kappa^2=.1$, which roughly corresponds to higher-frequency returns (say, daily returns).

As noted in CFS, the choice of the Gaussian auxiliary model, in conjunction with the constraints, means that  the GARCH(1,1) model is not well-equipped to handle the thick-tailed behavior exhibited by series generated from the log-normal SV model. Intuitively, this means that the constraints on the auxiliary parameters are likely to be binding since this auxiliary model is a crude approximation of the structural model. However, it is not certain if the inadequacy of the Gaussian GARCH(1,1) auxiliary model in this case, which was originally noted in Kim et al. (1998), is due to the model itself, the bindings constraints or a mixture of both issues. In this way, the FUNC based auxiliary estimator may be able to mitigate these issues since it captures, in some sense, the impact of the constraints.

The score based I-I objective function does not require a weighting matrix as we are in the just identified setting; i.e., we choose $W=I$. For computational simplicity, we fix the number of data replications to be $H=10$ across all Monte Carlo designs.\footnote{Optimization is carried out using an iterative Gauss-Seidel grid search approach. Starting values were obtained by first running a crude grid search over $\Theta$ and choosing the corresponding grid values that minimized the I-I objective function. Only one iteration of the minimization procedure was carried out and more efficient estimates could be obtained by considering multiple iterations.}  We illustrate the performance of our proposed I-I estimator across three different sample sizes, $T=500,1000,2000$, and consider $1000$ Monte Carlo replications for each sample size/parameter specification, leading to six separate specifications in total. 

\subsubsection{Monte Carlo Results}
\subsubsection*{Simulation Design one: $\theta^{0,1} = (-.736,.90,.363)^{\prime}$}
To understand the difference between the constrained and unconstrained auxiliary estimators, {Table \ref{Table1}} contains the frequency of binding constraints for the GARCH(1,1) auxiliary parameter estimates when we allow the boundary of the constrained space to drift, which replicates the behavior of a drifting DGP. Recalling that, in their assessment, CFS employ the constraint $\varphi\geq.025$, we choose a drifting boundary of $\varphi\geq \bar{\varphi}_{T}$ so that the drifting pseudo-true value $\beta^0_{T}$ lies in the interior of the constrained space, but $\beta^0$ can lie on the boundary. To ensure that our approach is comparable to the approach of CFS, we set $\bar{\varphi}_{T}:={T}^{-.5}$, which ensure that, at least for $T=2000$ the two constraints are comparable. 

For each replication, we calculate the auxiliary estimator $\hat{\beta}^{r}_{T}$ subject to the constraints in \eqref{GARCHineq}, where we require the constraint ${\varphi}\geq\bar{\varphi}_{T}$, and calculate $\widehat{\beta}_{T}$ by taking a Newton-step from $\hat{\beta}^{r}_{T}$. While no constraints are used in the calculation of  $\widehat{\beta}_{T}$ it is informative to ascertain the number of times this estimator would have caused the constraints to bind or be violated, as this will tell us, to some extent, what using the unconstrained $\widehat{\beta}_{T}$ buys us, at least in comparison with $\hat{\beta}^{r}_{T}$.

Table \ref{Table1} demonstrates that under the first design, the constraints for the auxiliary model are binding in a non-negligible portion on the replications. Interestingly, the FUNC estimator violates the constraint $\varphi\geq \bar{\varphi}_{T}$ less frequently than the constrained estimator. This is important given the results of Francq and Zakoian (2009), which demonstrate that the constraint $\varphi>0$ is misleading, suggesting to conclude ``that the GARCH(1,1) model is sufficient for financial data'', while additional ARCH lags would be relevant. As clearly explained by Francq and Zakoian (2009), this misleading conclusion is due to the fact that ``as a result of the positive constraints, it is possible that the fitted GARCH(1,1) models'' deliver a zero constrained estimator for the second ARCH lag while the score at this value is strongly positive. By construction, the FUNC estimator makes a Newton correction to take this positive score into account. Note that, in contrast, the FUNC estimator violates the stationarity constraint  $\varphi+\pi<1$ more  frequently than this constraint binds for the constrained estimator. It is worth realizing that this constraint may be irrelevant since, when the GARCH(1,1) model is misspecified, the fact  that the pseudo-true value violates the stationarity constraint does not imply that the process is non-stationary. This finding actually confirms a point made by Chib, Kim and Shephard (1998): when the data generating process is stochastic volatility (with true volatility persistence of $\delta^0=.90$), estimating the pseudo-true value for a GARCH(1,1) model can often deliver estimates of volatility persistence with $\varphi+\pi>.90$.

\begin{table}[htbp]
	\centering
	\caption{Binding constraints for auxiliary estimators $\hat{\beta}^{r}_{T}$ and $\widehat{\beta}^{}_{T}$ in design one: $\theta^{0,1}=(-.736,.90,.363)^{\prime}$. All terms are in percentages. For $\widehat{\beta}^{}_{T}$, the values represent the percentage where the FUNC estimator would have caused the constraint to bind or be violated.   }
	\begin{tabular}{rrrr|rrr|rrrr}
		\hline\hline
		&       &  T=500     &    &       &  T=1000  & & & T=2000& &\\
		\hline
		&       & $\hat{\beta}^{r}_{T}$ & $\widehat{\beta}^{}_{T}$ &       & $\hat{\beta}^{r}_{T}$ & $\widehat{\beta}^{}_{T}$ && $\hat{\beta}^{r}_{T}$ & $\widehat{\beta}^{}_{T}$ &\\\hline
		$\psi\geq0$   &       & 0.00\% & 1.20\% &       & 0.00\% & 0.20\% &&0.00\%&1.50\%&\\
		$\varphi\geq{T}^{-.5}$ &       & 11.10\% & 1.20\% &       & 4.50\% & 0.20\% &&0.60\%&0.20\%&\\
		$\pi\geq0$    &       & 3.80\% & 0.80\% &       & 0.40\% & 1.90\% &&0.70\%&1.10\%&\\
		$\varphi+\pi\leq1$ &       & 1.10\% & 1.90\% &       & 0.30\% & 0.50\% &&0.00\%&0.20\%&\\
		\hline\hline
	\end{tabular}%
	\label{Table1}%
\end{table}%

Summary statistics for the resulting score-based I-I parameter estimates are collected in Table \ref{Table2}. The results show that this I-I approach behaves well in finite samples, regardless of the constraints for the auxiliary model. To further understand the finite-sample properties of these estimators, in panels A and B of Figure \ref{Deltafig1} we plot the sampling distributions of the $\delta$ and $\sigma_{v}$ estimators across the three different sample sizes.\footnote{To ensure that all plots adequately represent the various sampling distributions and neatly fit in the same figure, we have thrown out 1.5\% of the lower tail observations for each series.} The results in Figure \ref{Deltafig1} are similar to those reported on page 963 in CFS.
\begin{table}[htbp]
	\centering
	\caption{Summary statistics for I-I estimates based on the proposed score approach in design one: $\theta^{0,1}=(-.736,.90,.363)^{\prime}$. STD- Monte Carlo standard deviation of the replications. RMSE- root mean squared error of the replications. M. Bias- mean bias of the replications.}
	\begin{tabular}{rrrr|rrr|rrr}
		\hline\hline
		T= 500    &       &       &       &   T=1000&             &       &    T=2000   &  &   \\
		\hline\hline
		$\theta$     &   STD   & RMSE  & M. Bias &               STD   & RMSE  & M. Bias & STD   & RMSE  & M. Bias\\\hline\hline
		$\alpha$ &  0.2849 & 0.2854 & 0.0178 &         0.1996 & 0.1996 & -0.0052 &0.1439&0.1439&-0.0011\\
		$\delta$  & 0.1299 & 0.1397 & -0.0503 &         0.0857 & 0.0897 & -0.0266 &0.0381&0.0392&-0.0092\\
		$\sigma_v$  & 0.0945 & 0.0961 & -0.0178 &         0.0447 & 0.0477 & -0.0003 &0.0333&0.0336&0.0047\\
		\hline\hline
	\end{tabular}%
	\label{Table2}%
\end{table}%

\subsubsection*{Simulation Design Two: $\theta^{0,2} = (-.141,.98,.0614)^{\prime}$}
Analyzing the frequency of binding constraints for the second Monte Carlo design, we find a very similar story to the first Monte Carlo design. Under this design, there are an even larger number of replications where the constraint $\varphi \geq \bar{\varphi}_{T}$ binds for the constrained estimator and is violated for the FUNC estimator.  As discussed by CFS, a small unconditional coefficient of variation for volatility creates a more challenging estimation problem, which seems to have had an impact on the frequency of binding constraints in this GARCH(1,1) auxiliary model. The aforementioned work of Francq and Zakoian (2009) may suggest that a GARCH(1,$q$), $q>1$, would have provided more informative parameters to match for I-I. Again, the fact that the FUNC estimator appears to be much less impacted by this constraint (five fewer times for a sample size of 500, and ten less times for larger sample sizes) is good news regarding its ability to capture the relevant information in the data. In addition, we find that there are a relatively large number of replications where the FUNC estimator does not satisfy the constraint $\varphi+\pi\leq1$: in about seven percent of the samples the FUNC estimator violates this constraint (at the sample size of $T=500$). This is not surprising since, as explained above, as discussed in Chib, Kim and Shephard (1998), we would expect the estimator of $\varphi+\pi$ to be larger than 0.98 (the true value of the volatility persistence, $\delta$, under this design).
\begin{table}[htbp]
	\centering
	\caption{Binding constraints for auxiliary estimators $\hat{\beta}^{r}_{T}$ and $\widehat{\beta}^{}_{T}$ in design two: $\theta^{0,2}=(-.141,.98,.0614)^{\prime}$. All terms are in percentages. For $\widehat{\beta}^{}_{T}$, the values represent the percentage where the FUNC estimator would have caused the constraint to bind or be violated.  }
	\begin{tabular}{rrrr|rrr|rrrr}
		\hline\hline
		&       &  T=500     &    &       &  T=1000  & & & T=2000& &\\
		\hline
		&       & $\hat{\beta}^{r}_{T}$ & $\widehat{\beta}^{}_{T}$ &       & $\hat{\beta}^{r}_{T}$ & $\widehat{\beta}^{}_{T}$ && $\hat{\beta}^{r}_{T}$ & $\widehat{\beta}^{}_{T}$ &\\\hline
		$\psi\geq0$   &       & 0.00\% & 6.10\% &       & 0.00\% & 2.40\% &&0.00\%&0.30\%&\\
		$\varphi\geq {T}^{-.5}$ &       & 29.50\% & 6.10\% &       & 20.70\% & 2.40\% && 10.90\%&0.30\%&\\
		$\pi\geq0$    &       & 0.00\% & 6.40\% &       & 0.00\% & 2.50\% &&0.00\%&0.50\%&\\
		$\varphi+\pi\leq1$ &       & 0.30\% & 6.50\% &       & 0.00\% & 2.60\% &&0.10\%&0.30\%&\\
		\hline\hline
	\end{tabular}%
	\label{Table3}%
\end{table}%

\begin{figure}[h]
	\centering
	\includegraphics[width=18cm, height=10cm]{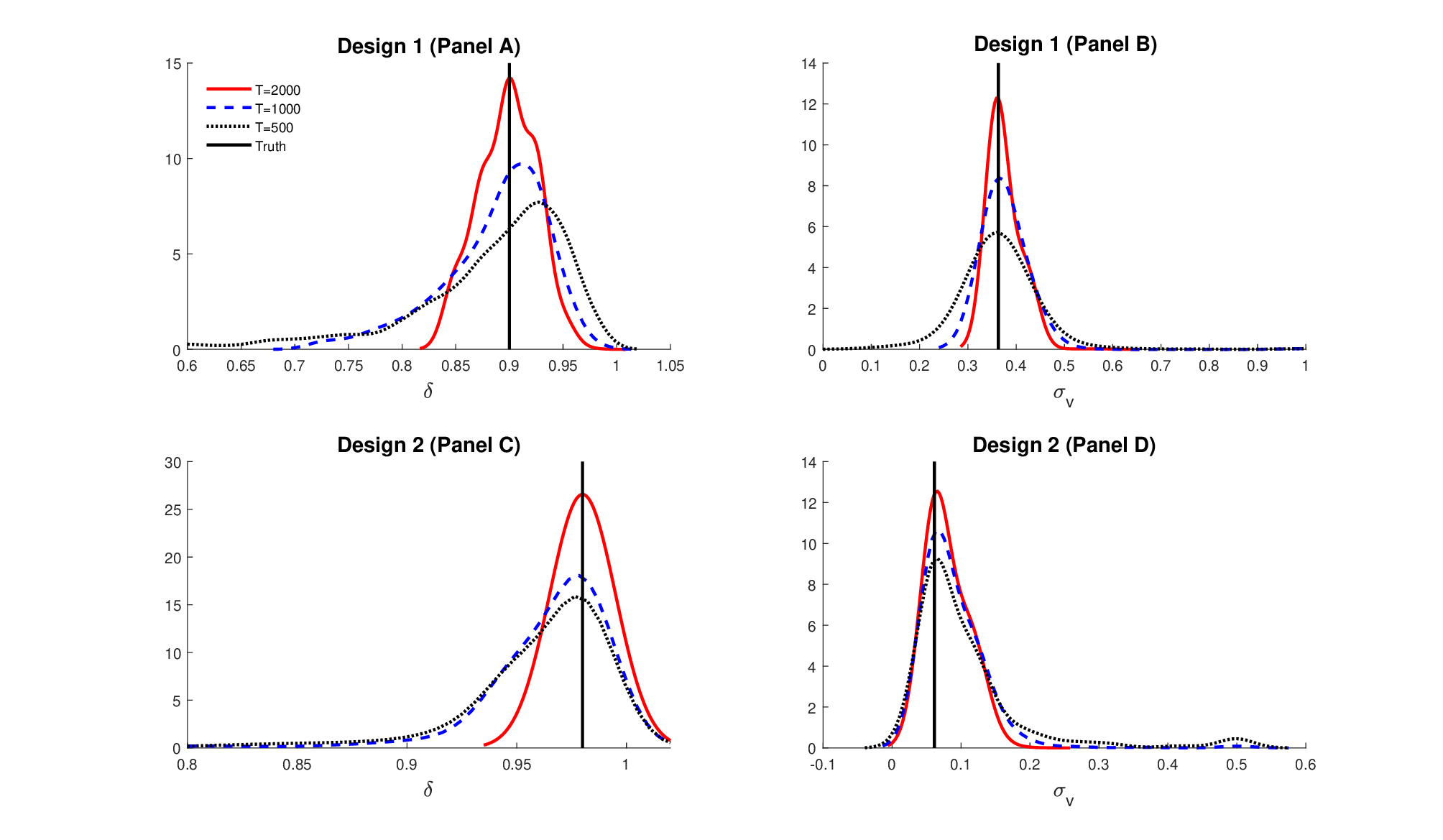}
	\caption{Sampling distribution for the I-I estimator of $\delta$ and $\sigma_v$ under Monte Carlo designs one, $\theta^{0,1}=(-.736,.90,.363)^{\prime}$, and two, $\theta^{0,2}=(-.141,.98,.0614)^{\prime}$. The line thick line corresponds to a sample size of T=2000, the thick dashed line corresponds to T=1000 and the thin dashed line corresponds to T=500. The vertical line represents the true value in the simulation.} 
	\label{Deltafig1}
\end{figure}

Summary statistics for the resulting score-based I-I parameter estimates are collected in Table \ref{Table4}, with the results reflecting the same conclusions as those obtained in the first Monte Carlo design. The sampling distributions of the $\delta$ and $\sigma_{v}$ estimators in the second Monte Carlo design are contained in panels C and D of Figure \ref{Deltafig1}. Again, the figures demonstrate that this approach works well.
\begin{table}[htbp]
	\centering
	\caption{Summary statistics for I-I estimates based on the proposed score approach in design two: $\theta^{0,2}=(-.141,.98,.0614)^{\prime}$. STD- Monte Carlo standard deviation of the replications. RMSE- root mean squared error of the replications. M. Bias- mean bias of the replications.}
	\begin{tabular}{rrrr|rrr|rrr}
		\hline\hline
		T= 500    &       &       &       &   T=1000&             &       &    T=2000   &  &   \\
		\hline\hline
		$\theta$     &   STD   & RMSE  & M. Bias &               STD   & RMSE  & M. Bias & STD   & RMSE  & M. Bias\\\hline\hline
		$\alpha$ &  0.5576 & 0.5576 & 0.0037 &         0.3857 & 0.3860 & -0.0171 &0.2822&0.2822&0.0018\\
		$\delta$  & 0.0592 & 0.0657 & -0.0287 &         0.0303 & 0.0321 & -0.0107 &0.0096&0.0102&-0.0034\\
		$\sigma_v$  & 0.0804 & 0.0805 & 0.0054 &         0.0358 & 0.0359 & 0.0030 &0.0097&0.0097&0.0008\\
		\hline\hline
	\end{tabular}%
	\label{Table4}%
\end{table}%

\subsection{$\alpha$-Stable Model}
In this section, we apply our I-I approach to data generated iid from the $\alpha$-stable distribution. We recall that the $\alpha$-stable distribution is characterized by four parameters:
$\alpha$- the tail index; $\gamma$- the skewness parameter; $\mu$- the location; and $\sigma$- the scale. To simplify the analysis, we assume that the location parameter $\mu$ is known and fix its value at $\mu=0$.\footnote{This is not overly restrictive as reliable estimators of the location parameter exist that can easily be employed before the analysis.}

The $\alpha$-stable distribution has no closed-form density representation and thus maximum likelihood estimation of the unknown parameters is difficult, which has led several authors to apply I-I to estimate the unknown parameters. Following the discussion in Section \ref{sec2}, we take as our auxiliary model the skewed Student-t (hereafter, skew-t) distribution of {Fernandez and Steel (1998)}: 
\begin{flalign*}
f(y;\beta)=\frac{\frac{\Gamma(\frac{\nu}{2}+\frac{1}{2})}{\Gamma(\nu/2)}}{\sqrt{\nu}}\frac{1}{\ell\left(\eta+\frac{1}{\eta}\right)}\left\{1+\frac{1}{\nu}\left(\frac{y-\omega}{\ell}\right)^2\left[\frac{1}{\eta^2}\1[y\geq \omega]+\eta^2\1[y<\omega]\right]\right\}^{-\frac{\nu+1}{2}},
\end{flalign*}where $\nu$ captures tail thickness, $\eta$ captures skewness, and the location and scale parameters are $\omega$ and $\ell$, respectively.

\subsubsection{Monte Carlo Design}

We consider a Monte Carlo design that is similar to Garcia et al. (2011). We fix $(\mu,\sigma,\gamma)'$ at $(0,0.5,0)'$ and consider two different true values for $\alpha^0$: 1) $\alpha^0=1.90$, and 2) $\alpha^0=1.95$. 
We consider two sample sizes, $T=500$ and $T=1000$, and generate 1000 Monte Carlo replications for each sample size.

As noted in Section \ref{sec2}, when $\alpha$ is close to, or larger than, 1.9, the parameter $\nu$ in the skew-t auxiliary model can become poorly identified in small samples and can result in I-I estimators with poor behavior.\footnote{A value of $\alpha<2$ implies that the unconditional variance of the random variable is not finite.} Simulation results in Garcia et al. (2011) demonstrate that, for data generated from the $\alpha$-stable distribution with $\alpha<2$ but close to $2$, the unconstrained PMLE for $\nu$ is non-Gaussian and numerically unstable in small samples. To circumvent this issue, the authors impose $\nu\leq 2$ within the skew-t auxiliary model, and then use as auxiliary statistics for I-I inference the auxiliary PMLE for $\beta$ and the corresponding KT multiplier $\lambda$ associated with the inequality constraint on $\nu$. {We refer the interested reader to Section two for further discussion on the need for this constraint.}

{Instead of imposing an arbitrary constraint on the auxiliary model, which may limit its identifying power, we propose to use the FUNC estimator.} In this context, the FUNC estimator does not display the same numerical instability observed in the unconstrained estimator for $\nu$: the calculation of the FUNC estimator is based on the numerically stable constrained estimator. Therefore, we argue that the FUNC estimator represents a practically useful medium between the fully unconstrained, and numerically unstable, auxiliary estimator and the standard constrained auxiliary estimator.

\subsubsection{Monte Carlo Results}
Under both Monte Carlo designs, the constraint on the auxiliary estimator ($\hat{\beta}^r_T$) for $\nu$ binds across all the replications, which implies that if one were to use the constrained I-I approach suggested by CFS, identification of $\alpha^0$ is completely determined by the KT multiplier associated with the constraint. Likewise, across all replications and both Monte Carlo designs, the FUNC estimator ($\widehat{\beta}_{T}$) would have caused the constraints to bind or be violated for both Monte Carlo designs estimator. Therefore, we can conclude that the constrained version of this auxiliary model is unable to identify $\alpha^0$ by itself, and instead must rely on additional auxiliary statistics that are not associated with the auxiliary {model; i.e., the KT multiplier in this case.}

In contrast, the FUNC estimator only uses the information contained in the (unconstrained) auxiliary model. Since the FUNC estimator is unconstrained, it is not surprising that this estimated auxiliary parameter always violates the constraint $\nu\leq2$. Given the above results for the constrained and FUNC-based auxiliary estimators, it is clear that the pseudo-true value defined by this constrained optimization program is on the boundary of the parameter space. As such, an alternative strategy would be to consider a skew-t auxiliary model that imposes the equality constraint $\nu=2$.

{Before presenting the Monte Carlo results, we remark that if the pseudo-true value of the auxiliary parameters were to violate the constraint $\nu\leq2$, the auxiliary score (evaluated at this pseudo-true value) would not be zero in the limit, and, as such, this case would be outside the scope of our theoretical analysis. While such scenarios could be accommodated by extending our theoretical framework, at the cost of additional notation and technical arguments, this extension is not germane to the main message of the paper. Therefore, for the sake of brevity, we do not consider such cases further.} 

Table \ref{stable_table2} reports the standard deviation (STD), root-mean-squared error (RMSE) and mean bias (M. Bias) associated with the parameter estimates from our unconstrained I-I estimation approach; namely, $\alpha,\gamma,\sigma$. The results demonstrate that our approach produces estimates with reliable finite-sample properties across both sample sizes. 
\begin{table}[htbp]
	\centering
	\caption{Summary statistics for Monte Carlo estimators across both sample sizes for design one ($\alpha^0=1.90$) and two ($\alpha^0=1.95$), denoted as D1 and D2 in the table. STD- Monte Carlo standard deviation of the replications. RMSE- root mean squared error of the replications. M. Bias- mean bias of the replications. Across both designs $\gamma^0=0$ and $\sigma^0=.5$. Under both designs, the constrain $\nu\leq2$ was binding in all Monte Carlo replications for both $\hat{\beta}^r_T$ and $\widehat{\beta}_{T}$. }
	\begin{tabular}{rrrr|rrr}
		\hline\hline
		&&         T=500     &         &    & T=1000 &\\\hline\hline
		D1 &  $\alpha$      & $\gamma$ &   $\sigma$  &  $\alpha$      & $\gamma$&$\sigma$\\\hline\hline
		STD    &0.0817& 0.0003&0.0295 &0.0612 &0.0002 &0.0210\\
		RMSE   &0.0819& 0.0003&0.0295 &0.0611 &0.0003 &0.0210\\
		M.Bias &-0.0060& 0.0002&-0.0016 &-0.0016 &0.0001 &0.0007\\\hline\hline
		D2	&  $\alpha$      & $\gamma$ &   $\sigma$  &  $\alpha$      & $\gamma$&$\sigma$\\\hline\hline
		STD &0.0637& 0.0002&0.0272 &0.0471 &0.0002 &0.0195\\
		RMSE &0.0646& 0.0003&0.0272 &0.0474 &0.0003 &0.0195\\
		M.Bias &-0.0108& 0.0002&-0.0007 &-0.0058 &0.0001 &0.0005\\\hline\hline
	\end{tabular}%
	\label{stable_table2}%
\end{table}%

In addition, Figure \ref{stable_fig1} contains kernel density estimates of the standardized parameter estimates across the Monte Carlo replications for the case of $\alpha^0=1.95$, the corresponding results for the design where $\alpha^0=1.90$ are very similar and are not reported for the sake of brevity. The results demonstrate that the standardized estimators have a roughly Gaussian shape, even though they are computed from $\alpha$-stable random variables. We refer the reader to Garcia et al. (2011) for theoretical justification of this phenomena. 
\begin{figure}[h!]
	\centering 
	\setlength\figureheight{3.65cm} 
	\setlength\figurewidth{5.75cm} 
	\input{Alpha_Stab.tikz} 
	\caption{Sampling distribution for the standardized I-I estimator of $\theta$ in the $\alpha$-stable Monte Carlo experiments: $\theta^{0}=(\alpha^0,\gamma^0,\sigma^0)'=(1.95,.0,.5)^{\prime}$. The results under the Monte Carlo experiment with $\theta^0=(1.90,.0,..5)^{\prime}$ are similar, and therefore not presented for brevity.} 
	\label{stable_fig1} 
\end{figure}

\subsection{Stochastic Volatility Jump Diffusion (SVJD) Model}
Motivated by the now well-established empirical findings of time-varying volatility and the existence of jumps in returns data, we explore here a continuous-time specification for financial returns. We consider that returns evolve in continuous-time according to a mean reverting stochastic volatility model, which follows an Ornstein-Uhlenbeck process, and where returns themselves exhibit random jumps. We first demonstrate, via Monte Carlo results, that in this empirically relevant model, the class of GARCH auxiliary models will often deliver estimated auxiliary parameters that are near the boundary of the parameter space. We then demonstrate that our I-I approach delivers reliable estimators of the corresponding structural parameters even though the original auxiliary parameters are near the boundary of the parameter space. Lastly, we use our I-I approach to conduct inference on the parameters of the SVJD model to determine whether or not there exists significant evidence of jumps in daily S\&P500 returns.

For $P_{t}$ denoting the asset price at time $t>0$, and $p_{t}:=\ln(P_{t})$, assume that $p_t$ evolves according to the bivariate diffusion process
\begin{flalign} d p _ { t } & = { \mu } d t + \exp( { V _ { t } }/2)  d W _ { t } ^ { p } + d J _ { t } ^ { p }, \nonumber\\ d V _ { t } & =  { \kappa } \left(  { \eta } - V _ { t } \right) d t +  { \sigma } _ { v } d W _ { t } ^ { v },\nonumber\\dJ_{t}^{p}&=Z_{t}d N_{t},\;\; Z_{t}\sim N(\mu_{j},\sigma^2_{j}),\;\;\text{Pr}\left[d N_{t}=1\right]=\bar{\lambda}_{j}dt+o(t)\label{svj_simple},
\end{flalign}where $dW _ { t } ^ { v }$ and $dW _ { t } ^ { p }$ are correlated Brownian motion processes, with correlation $\rho$, and $dN_{t}$ is a Poisson process with intensity $\bar{\lambda}_{j}$, and $Z_{t}\sim N(\mu_{j},\sigma^2_{j})$. We collect the unknown parameters into $\theta=(\mu,\kappa,\eta,\sigma^2_{v},\bar{\lambda}_{j},\mu_{j},\sigma^2_{j},\rho)'$ and consider inference on $\theta$ using I-I. 

Following the analysis in Example 1, we take as our auxiliary model the GARCH model with Student-t errors: for $r_t:=p_t-p_{t-1}$ denoting log-returns,
\begin{flalign*}
r_{t}&=\mu_a+\sqrt{h_{t}}\epsilon_{t},\\h_{t}&=\psi+\varphi \left(r_{t-1}-\mu_a\right)^2+\pi h_{t-1},
\end{flalign*}{where, for $v(\eta):=[{{(1/\eta-1/2)}/{1/\eta}}]$, $\epsilon_{t}\sim_{iid} v(\eta)^{1/2}t_{1/\eta}$, and $v(\eta)^{1/2}t_{1/\eta}$ denotes a Student-t with unit variance.}  The auxiliary GARCH model is unable to identify the jump parameters in \eqref{svj_simple}, and so we supplement the auxiliary GARCH model
with additional summary statistics based on both bipower variation and realized jump variation: for $r_{t,{i}}$ denoting the $i$-th, out of $M$, equally-spaced intra-day returns  observed on day $t$, bipower variation is defined as $$B V _ { t } := \frac { \pi } { 2 } \left( \frac { M } { M - 1 } \right) \sum _ { i = 2 } ^ { M } \left| r _ { t ,{ i } } \right| \left| r _ { t ,{( i - 1) } } \right|,$$ and jump-variation is defined as $$\text{JV}_{t}:=\max\{\text{RV}_{t}-\text{BV}_{t},0\},\;\;\text{RV}_{t}=\sum_{i=1}^{M}r_{t,{i}}^2,$$ where $\text{RV}_{t}$ denotes realized volatility. For I-I estimation we then consider the additional summary statistics: for $\overline{\text{JV}}:=\frac{1}{T}\sum_{t=1}^{T}\text{JV}_{t}$, 
\begin{flalign*}
S_{1}&:=\frac{1}{T}\sum_{t=1}^{T}\text{sign}(r_{t})\sqrt{\text{JV}_{t}},\;
S_{2}:=\frac{1}{T}\sum_{t=1}^{T}\left(\text{JV}_{t}-\overline{\text{JV}}_t\right)^2,\;
S_{3}:=\frac{1}{T}\sum_{t=2}^{T}\left(\text{JV}_{t}-\overline{\text{JV}}_t\right)\left(\text{JV}_{t-1}-\overline{\text{JV}}_t\right)
\end{flalign*}which correspond to the mean, variance and covariance of the realized jump variation. We note here that {Frazier et al. (2018)} have used these statistics to help identify the jump-process parameters in a discrete-time version of this SVJD model for daily S\&P500 returns data.

Estimating GARCH models on daily data can often lead to estimated values of $\varphi$ and $\pi$ such that the constraint $\varphi+\pi\leq1$ is very close to binding. To ensure numerical stability of the optimization procedure, in practice this constraint is often implemented as $\varphi+\pi\leq1-c$, for $c>0$ and small. 

With this point in mind, and similar to the $\alpha$-stable example, there is every reason to suspect that a lack of variation in the estimated constrained auxiliary parameters of $\varphi$ and $\pi$, due to the need to satisfy the constraint $\varphi+\pi\leq1$, may cause identification issues for the I-I estimator of the structural parameters; i.e., if there are many parameter combinations for $\theta$ that yield estimated constrained auxiliary parameters close to the boundary of the constraints, then those simulated data sets may be difficult for the I-I objective function to distinguish. Therefore, we argue that an appropriate approach to the use of GARCH models in such settings is to employ our FUNC-based II approach. 

\subsubsection{Monte Carlo Design}
The observed data is generated from the model in equation \eqref{svj_simple} using an Euler discretization scheme with step size $1/\delta$. Log-prices are then generated according to the recursive scheme
\begin{flalign*}
p_{t,{(i+1)}/\delta}&=p_{t,{i}/\delta}+\mu\frac{1}{\delta}+\exp(V_{t,i/\delta}/2)\epsilon^{p}_{t,i}\frac{1}{\delta}+Z_{t,i}\Delta N_{t,i/\delta},\\
V_{t,{(i+1)}/\delta}&=V_{t,{i}/\delta}+\kappa(\eta-V_{t,i/\delta})\frac{1}{\delta}+\frac{\sigma_{v}}{\delta}\left(\rho\epsilon^{p}_{t,i}+\sqrt{1-\rho^2}\epsilon^{v}_{t,i}\right),
\end{flalign*}where $(\epsilon^{p}_{t,i},\epsilon^{v}_{t,i})'$ is bivariate standard normal, $Z_{t,i}\sim_{iid}N(\mu_j,\sigma^2_j)$ and $\Delta N_{t,i/\delta}$ is drawn from a Poisson distribution with intensity $\bar{\lambda}_{j}/\delta$. 

We simulate data from the above Euler approximation at the (approximate) one minute frequency, $\delta=400$, and record daily returns, and ten minute intra-day returns. The remaining simulations are discarded. We retain $T=500$ trading days for the Monte Carlo, which amounts to approximately two years of daily returns. For the intra-day returns, we use ten-minute returns to calculate our measures of realized volatility $\text{RV}^{10}_{t}:=\sum_{t=1}^{M}r_{t,{i}}^2$, $M=40$, and our corresponding realized jump measure $\text{JV}_{t}$. For the data simulation, we consider an initialization period of 400 periods, or one trading day. Following the design of {Creel and Kristensen (2015)}, we set the unknown parameters to $\theta^0=(0.00,0.02,0.25,0.20,0.10,0.00,0.50,-0.10)'$. 

For I-I estimation, we also consider a Euler discretization scheme at the one-minute sampling frequency, although finer-sampling schemes may lead to estimators with better properties. Since this scheme already requires simulating a large number of data points, to ease the computational burden of the I-I estimator, we limit the analysis to consider only a single simulated path, i.e., $H=1$. This will induce some efficiency loss in the resulting estimators, however, that is the price to pay for computational convenience.

\subsubsection{Monte Carlo Results}
First, we analyze the frequency at which the estimated auxiliary parameters cause the constraint $\varphi+\pi\leq1$ to bind.  From Table \ref{SVJD_tab} we see that in 26\% of the simulations the constraint was binding for the constrained auxiliary estimators, while the FUNC estimator either led to a value of the estimated parameters that was greater than or equal to one. As explained in Section \ref{sv_ex}, a true daily volatility persistence $\exp(-\kappa) \approx .98$ leads us to expect an estimator of $\varphi+\pi$ even larger than .98. Moreover, this behavior is amplified in the presence of jumps, and  leads to an unconstrained estimations of  $\varphi+\pi$ that if frequently larger than unity.

While this behavior of the FUNC estimator may seem counter-intuitive, recall that the condition $\varphi+\pi<1$ is required for stationarity in the GARCH model, under the explicit assumption that the true DGP is GARCH. Therefore, in this example the constraint $\varphi+\pi<1$ is meaningless since the true DGP is not GARCH. Indeed, the process in this simulation is stationary even though the estimated values can satisfy $\varphi+\pi\geq1$ in any finite-sample. This example clearly demonstrates that if we simply used the constrained GARCH auxiliary estimators, we could artificially limit the identifying power of this auxiliary model. 

Using our FUNC-based I-I approach, we estimate the structural parameters of the SVJD model and report the results in Table \ref{SVJD_tab}. The corresponding estimators display low mean bias and reasonable values for the RMSE, especially given that we have used only a single simulated path for I-I. Similar to the previous two examples, the results of this section demonstrate that our approach is able to achieve identification without the need to resort to an I-I approach that utilizes the KT multipliers associated with the constraints for the auxiliary model.

\begin{table}[htbp]
	\centering
	\caption{Summary statistics for Monte Carlo estimators for the SVJD example. The sample size is $T=500$. STD- Monte Carlo standard deviation of the replications. RMSE- root mean squared error of the replications. M. Bias- mean bias of the replications. For this Monte Carlo example, the constraint on $\varphi+\pi\leq 1$ was binding in 26\% of the Monte Carlo replications, for both the constrained estimator $\hat{\beta}^r_{T}$ and the FUNC estimator $\widehat{\beta}_{T}$. }
	\begin{tabular}{rrrrrrrrr}
		\hline\hline
		$\theta$  &  $\mu$      & $\kappa$ &   $\eta$  &  $\sigma_v$      & $\bar{\lambda}_j$&$\mu_{j}$&$\sigma^2_j$&$\rho$\\\hline\hline
		$\theta^0$		& 0.0000  &  0.0200 &   0.2500 &   0.2000 &         0.1000 &  0.0000& 0.5000  & -0.1000\\ \hline\hline
		STD &0.0010 &   0.0068 &   0.1005 &   0.1193 &      0.0195& 0.0010 &  0.1099&    0.0204\\
		RMSE &0.0010  &  0.0068 &   0.1006 &   0.1193 &      0.0197&  0.0010 &  0.1098 &   0.0204\\
		M.Bias &0.0000 &  -0.0002 &   0.0062 &  -0.0048  &    0.0027 & 0.0000 &-0.0005 &  -0.0014
		\\\hline\hline
		
	\end{tabular}%
	\label{SVJD_tab}%
\end{table}%

\subsubsection{Empirical Illustration: S\&P500 Data} 
To further illustrate our approach, we apply our method to the SVJD model with leverage effects based on demeaned S\&P500 returns observed at the daily frequency, between 3 January 2017 and 3 January 2019, which consists of 501 daily observations. High-frequency intra-day returns are used to build realized volatility and bi-power variation estimators at the five minute sampling frequency. We source the data from the Oxford-Man Institutes ``realised library'', which contains daily returns on several important financial indices, and non-parametric volatility measures (Gerd et al., 2009).

The auxiliary model is again taken to be the GARCH model with standardized Student-t errors (and unit variance): for $r_t:=p_t-p_{t-1}$ denoting log-returns on the S\&P500 index, with $t=1,\dots,501$ denoting the daily frequency,
\begin{flalign*}
r_{t}&=\mu_a+\sqrt{h_{t}}\epsilon_{t},\\h_{t}&=\psi+\varphi \left(r_{t-1}-\mu_a\right)^2+\pi h_{t-1}.
\end{flalign*} Again, for $\overline{\text{JV}}:=\frac{1}{T}\sum_{t=1}^{T}\text{JV}_{t}$, we augment this model with the summary statistics 
\begin{flalign*}
S_{1}&:=\frac{1}{T}\sum_{t=1}^{T}\text{sign}(r_{t})\sqrt{\text{JV}_{t}},\;
S_{2}:=\frac{1}{T}\sum_{t=1}^{T}\left(\text{JV}_{t}-\overline{\text{JV}}\right)^2,\;
S_{3}:=\frac{1}{T}\sum_{t=2}^{T}\left(\text{JV}_{t}-\overline{\text{JV}}\right)\left(\text{JV}_{t-1}-\overline{\text{JV}}\right)
\end{flalign*}that capture the jump components of the SVJD model. While this simple SVJD model has now been generalized in several directions, e.g., with the inclusion of autocorrelated jumps (Fulop et al., 2014, Ait-Sahalia et al., 2015, and Maneesoonthorn et al., 2017), this simpler SVJD model is still empirically relevant since, if the corresponding jumps components are not statistically significant, it is highly unlikely that these more complicated modeling approaches are necessary.

Maximum likelihood-based inference on the full set of static parameters in this SVJD model, $$\theta=(\kappa,\eta,\sigma^2_{v},\bar{\lambda}_{j},\mu_{j},\sigma^2_{j},\rho)',$$
is challenging, due to the existence of the latent volatilities.\footnote{Since the data is demeaned before hand, we do not estimate $\mu$ in this example.} In contrast, I-I is straightforward due to our ability to cheaply simulate data from this model. Therefore, we consider estimation of the model in \eqref{svj_simple} using our I-I approach base on the GARCH auxiliary model, and where the statistics $(S_1,S_2,S_3)'$ yield auxiliary moments that enable us to identify the jump components. {Similar to the Monte Carlo example, the I-I procedure uses a simulation frequency of one minute.} 

The estimated values of $\theta$ obtained using this setup and the daily S\&P500 data are given in Table \ref{SVJD_tab_emp}, along with the corresponding standard errors. The standard errors are calculated using a block bootstrap approach, with 999 bootstrap replications and with a block length of twenty-five observations. Given the relatively short length of the time series, we believe these bootstrap standard errors are likely more reliable than those obtained from the asymptotic formula presented in Section three. 

Before analyzing the results in the table, we first note that the estimated auxiliary parameters for the GARCH model are such that $\varphi+\pi\approx1$ for both the constrained and FUNC-based auxiliary estimators. {Recall that, even though in this example the FUNC and constrained auxiliary estimators are similar, the FUNC estimator is guaranteed to be asymptotically normal, whereas the constrained estimator will in general not be asymptotically normal. Therefore, we contend that, even though the two estimators are similar in this small scale example, the use of the FUNC estimator for I-I is a safer choice than the constrained estimator.} 

Analyzing the results for $\theta$ in Table \ref{SVJD_tab_emp}, we see that the majority of the coefficients are statistically significant and have the correct signs and magnitudes, with the majority of the results being similar to those obtained elsewhere; see, e.g., Creel and Kristensen (2015). In particular, the results suggest that the jump process has a significant jump frequency, but that the resulting jump sizes are small, negative, and have large variability. {In addition, these results also suggest that the resulting jump components $(\mu_j,\sigma^2_j)$ are difficult to accurately measure.}

\begin{table}[htbp]
	\centering
	\caption{Estimators (Est) and standard errors (STD) for the SVJD S\&P500 exercise. For this data set, the constraint on $\varphi+\pi\leq 1$ was binding for both the constrained estimator $\hat{\beta}^r_{T}$ and the FUNC estimator $\widehat{\beta}_{T}$. }
	\begin{tabular}{rrrrrrrr}
		\hline\hline
		$\theta$  &   $\kappa$ &   $\eta$  &  $\sigma_v$      & $\bar{\lambda}_j$&$\mu_{j}$&$\sigma^2_j$&$\rho$\\\hline\hline
		Est & 0.2876  &  0.5945 &   0.1166  &  0.1278 &  -0.0034  &  1.2650 &  -0.6131\\
		STD & 0.0836  &  0.0539  &  0.0366 &  0.0383  &   0.0203  &   0.2902 &  0.1717\\
	\hline\hline
		
	\end{tabular}%
	\label{SVJD_tab_emp}%
\end{table}%

\section{Conclusion}

The overall message of this paper can be summarized as follows:
application of the I-I methodology may require the imposition of certain constraints on the auxiliary parameters, however, one must bear in mind that the behavior of I-I estimators for the structural parameters can be adversely affected by the constraints placed on the auxiliary parameters. In place of these constrained auxiliary parameters, our proposed strategy is {to use, for the purpose of I-I,} a FUNC (Feasible
UNConstrained) estimator of the auxiliary parameters, which, in spite of
being unconstrained, is always well-defined. 

This FUNC estimator leads to simple score and Wald-based I-I approaches, which have been shown to be asymptotically equivalent, at least to first-order, with the approach based on constrained auxiliary parameters proposed by CFS. Several Monte Carlo studies demonstrate the good finite-sample properties of this approach, and document that our I-I estimator can deliver robust estimators of the corresponding structural parameters, even in cases where the pseudo-true value of the auxiliary parameters is on the boundary of the parameter space.

\appendix 
\section{Proofs of Main Results}
\begin{proof}[Proof of Lemma 1]
We first prove that $\|\hat{\beta}^r_{T}-\beta^0_{T}\|=o_{P}(1)$. The argument follows the standard approach. Under $\{\theta_{T}\}\in\Gamma(\theta^0,\beta^0)$,
\begin{eqnarray*}
0&\leq&\mathcal{Q}(\theta_{T},{\beta}^0_{T})-\mathcal{Q}(\theta_{T},\hat{\beta}^{r}_{T}) \\
&=&\mathcal{Q}(\theta_{T},{\beta}^0_{T})-\mathcal{Q}(\theta_{T},\hat{\beta}^{r}_{T})+Q_{T}(\hat{\beta}^{r}_{T})-Q_{T}(\hat{\beta}^{r}_{T})+Q_{T}({\beta}^0_{T})-Q_{T}({\beta}^0_{T})\\ &\leq& 2\sup_{\beta\in\mathbf{B}}|Q_{T}(\beta)-\mathcal{Q}(\theta_{T},\beta)|+o_{P}(1),
\end{eqnarray*}where the $o_{P}(1)$ follows from the fact that $Q_{T}(\hat{\beta}^{r}_{T})\geq Q_{T}(\beta^0_{T})+ o_{P}(1)$. From the uniform convergence in \textbf{A0}(i), we can conclude 
$$0\leq \mathcal{Q}(\theta_{T},{\beta}^0_{T})-\mathcal{Q}(\theta_{T},\hat{\beta}^{r}_{T})\leq o_{P}(1).$$The result then follows from the identification condition in \textbf{A0}(ii).

We next show $\sqrt{T}\left(\hat{\beta}^r_{T}-\beta^0_{T}\right)=O_{P}(1).$ The quadratic expansion 
	\begin{equation*}
	Q_{T}(\beta )=Q_{T}(\beta _{T}^{0})+\frac{\partial Q_{T}(\beta _{T}^{0})}{%
		\partial \beta ^{\prime }}(\beta -\beta _{T}^{0})+\frac{1}{2}(\beta -\beta
	_{T}^{0})^{\prime }\frac{\partial ^{2}Q_{T}(\beta _{T}^{0})}{\partial \beta
		\partial \beta ^{\prime }}(\beta -\beta _{T}^{0})+R_{T}(\beta ), 
	\end{equation*}can be rewritten as 
	\begin{equation}\label{new_exp}
	Q_{T}(\hat{\beta}^r_{T} )=Q_{T}(\beta _{T}^{0})+\frac{1}{T}\kappa_{T}'J_{T}^{-1/2}\sqrt{T}\frac{\partial Q_{T}(\beta^0_T)}{\partial\beta}-\frac{1}{2T}\|\kappa_{T}\|^2+R_{T}(\hat{\beta}^{r}_{T} ), 
	\end{equation}	 where $$J_{T}=-\frac{\partial^2 Q_{T}(\beta^0_{T})}{\partial\beta\partial\beta'},\; \kappa_{T}=J_{T}^{1/2}\sqrt{T}(\hat{\beta}^{r}_{T}-\beta^0_{T}).$$From the definition of $\hat{\beta}_{T}^{r}$, applying the quadratic expansion in \eqref{new_exp}
	\begin{eqnarray*} o _ { P } ( 1 ) & \leq& T\cdot\left[Q _ { T } ( \hat { \beta }^{r}_{T} ) - Q_ { T } \left( \beta_{T}^ { 0 } \right) \right]\\ & = &\kappa _ { T } ^ { ' } J _ { T } ^ { -1 / 2 } \sqrt{T}\frac{\partial Q_{T}(\beta^0_{T})}{\partial\beta} - \frac { 1 } { 2 } \left\| \kappa _ { T } \right\| ^ { 2 } + T\cdot R _ { T } ( \hat { \beta }^{r}_{T} ) 
		\end{eqnarray*}
However, by Assumption \textbf{A2}, since $\|	\hat { \beta }^{r}_{T}-{ \beta }^{0}_{T}\|=o_{P}(1)$, we have that $$\left| T\cdot R _ { T } ( \hat { \beta }^{r}_{T} )\right|\leq \left(1+\left\|\sqrt{T}\left(\hat{\beta}^r_{T}-\beta^0_T\right)\right\|\right)^2o_{P}(1).$$ Applying the above, we have
\begin{eqnarray*}
&&\kappa _ { T } ^ { ' } J _ { T } ^ { -1 / 2 } \sqrt{T}\frac{\partial Q_{T}(\beta^0_{T})}{\partial\beta} - \frac { 1 } { 2 } \left\| \kappa _ { T } \right\| ^ { 2 } + T\cdot R _ { T } ( \hat { \beta }^{r}_{T} ) \\	& =& O _ { P } \left( \left\| \kappa _ { T } \right\| \right) - \frac { 1 } { 2 } \left\| \kappa _ { T } \right\| ^ { 2 } + \left( 1 + \left\| J _ { T } ^ { - 1 / 2 } \kappa _ { T } \right\| \right) ^ { 2 } o _ { P } ( 1 ) \\ & = &O _ { P } \left( \left\| \kappa _ { T } \right\| \right) - \frac { 1 } { 2 } \left\| \kappa _ { T } \right\| ^ { 2 } + o _ { P } \left( \left\| \kappa _ { T } \right\| \right) + o _ { P} \left( \left\| \kappa _ { T } \right\| ^ { 2 } \right) + o _ { P } ( 1 ) ,
\end{eqnarray*}
so that we may rewrite the above as
	$$\|\kappa_{T}\|^2\leq2 \|\kappa_{T}\|\left[o_{P}(1)+O_{P}(1)\right]+o_{P}(1).$$ Defining $\zeta_{T}:=\left[o_{P}(1)+O_{P}(1)\right]\equiv O_{P}(1)$ and $x_{T}:=\|\kappa_{T}\|$, we end up with the inequality
	\begin{eqnarray*}
	x_{T}^2-2x_{T}\zeta_{T}+o_{P}(1)\leq0, 	
	\end{eqnarray*}which is satisfied for $x_{T}$ in the interval $\zeta_{T}\pm\sqrt{\zeta^2_T+o_{P}(1)}$; i.e., for $x_{T}$ in the interval $[0,2\zeta_{T}+o_{P}(1)]$. Hence, $\zeta_T=O_{P}(1)$ implies that $$x_{T}=\|\kappa_{T}\|=\left\|J_{T}^{1/2}\sqrt{T}\left(\hat{\beta}^r_{T}-\beta^0_T\right)\right\|=O_{P}(1)$$ which implies that
$$\sqrt{T}\left(\hat{\beta}^r_{T}-\beta^0_T\right)=O_{P}(1).$$

	Now, we prove $\sqrt{T}\hat{\lambda}_{T}=O_{P}(1)$. First, consider the case where $q>d_\beta$. By assumption, at most there are $\tilde{q}$ dimensions of $g({\beta}^0_T)$ that are precisely zero, which are all contained in the vector $\tilde{g}(\beta^0_T)$. As shown above, $\sqrt{T}(\hat{\beta}^r_T-\beta^0_T)=O_{p}(1)$, and we can then be sure that asymptotically, with probability one, all zero entries of $g(\hat{\beta}^r_T)$ are also included in $\tilde{g}(\hat{\beta}^r_T)$.  Define $\tilde{\lambda}_T$ to be the $\tilde{q}$-dimensional sub-vector of $\hat{\lambda}_T$ that corresponds to the entries of $g$ that are in $\tilde{g}$. By the slackness conditions of the Kuhn-Tucker optimization problem, asymptotically, with probability one, since  $\tilde{g}(\hat{\beta}^r_T)$ contains all the zero entries of ${g}(\hat{\beta}^r_T)$, $\tilde{\lambda}_T$ contains all the possible non-zero entries of $\hat{\lambda}_T$.

	With these definitions, the Kuhn-Tucker first-order conditions can be stated as
	\begin{eqnarray*}
	\sqrt{T}\frac{\partial Q_{T}(\hat{\beta}^{r} _{T})}{\partial \beta ^{\prime }}+\frac{\partial \tilde{g}^{\prime }(\hat{\beta}^{r} _{T})}{\partial \beta }%
		\sqrt{T}\tilde{\lambda}_{T}&=&0.
	\end{eqnarray*}
For some intermediate value $\bar{\beta}_T$, a first-order expansions gives
\begin{eqnarray*}
\sqrt{T}\frac{\partial Q_{T}(\beta _{T}^{0})}{\partial \beta ^{\prime }}+%
\frac{\partial ^{2}Q_{T}\left(
	\bar{\beta} _{T}^{}\right)}{\partial \beta \partial \beta ^{\prime }}\sqrt{T} \left( \hat{\beta}_{T}^{r}-\beta
_{T}^{0}\right) +\frac{\partial \tilde{g}^{\prime }({\beta} _{T}^{0})}{\partial \beta }%
\sqrt{T}\tilde{\lambda}_{T}+o_P(\sqrt{T}\tilde{\lambda}_T) &=&0,
\end{eqnarray*}where the $o_P(\sqrt{T}\tilde{\lambda}_T)$ term follows by the first part of Lemma 1 and the continuity of $\partial g(\beta)/\partial\beta'$ in \textbf{A1}(iv). Since ${\partial \tilde{g}^{\prime }({\beta}_{T}^{0})}/{\partial \beta }$ is
full column-rank, we have, for some intermediate value $\bar{\beta}_{T}$,
\begin{eqnarray*}
\sqrt{T}\tilde{\lambda}_{T} &=&-\left[ \frac{\partial \tilde{g}({\beta}_{T}^{0})}{%
	\partial \beta ^{\prime }}\frac{\partial \tilde{g}^{\prime }({\beta}_{T}^{0})}{%
	\partial \beta }\right] ^{-1}\frac{\partial \tilde{g}({\beta}_{T}^{0})}{\partial
	\beta ^{\prime }}\sqrt{T}\frac{\partial Q_{T}(\hat{\beta}_{T}^{r})}{\partial
	\beta }+o_P(\sqrt{T}\tilde{\lambda}_T) \\
&=&-\left[ \frac{\partial \tilde{g}({\beta}_{T}^{0})}{\partial \beta ^{\prime }}%
\frac{\partial\tilde{g}^{\prime }({\beta}_{T}^{0})}{\partial \beta }\right]
^{-1}\frac{\partial \tilde{g}({\beta}_{T}^{0})}{\partial \beta ^{\prime }}\left[ 
\sqrt{T}\frac{\partial Q_{T}(\beta _{T}^{0})}{\partial \beta }+\frac{%
	\partial ^{2}Q_{T}(\bar{\beta}_{T})}{\partial \beta \partial \beta ^{\prime }%
}\sqrt{T}\left( \hat{\beta}_{T}^{r}-\beta _{T}^{0}\right) \right] +o_P(\sqrt{T}\tilde{\lambda}_T)\\&=&O_{P}(1)+o_P(\sqrt{T}\tilde{\lambda}_T),
\end{eqnarray*}where the last line follows from Assumption \textbf{A1}(iii) and $\sqrt{T}(\hat{\beta}^r_T-\beta^0_T)=O_{P}(1).$

In the case where $q\leq d_\beta$, the above arguments remain valid if we take $\tilde{q}=q$, $\tilde{g}(\beta)=g(\beta)$, $\tilde{\lambda}_T=\hat{\lambda}_T$, and note that ${\partial \tilde{g}^{\prime }({\beta}_{T}^{0})}/{\partial \beta }$ has full column-rank $q$. 
\end{proof}
	
\begin{proof}[Proof of Proposition 1]
 A first-order expansions of the first-order conditions \eqref{lagrange} gives
\begin{eqnarray}
\sqrt{T}\frac{\partial Q_{T}(\beta _{T}^{0})}{\partial \beta ^{\prime }}+%
\frac{\partial ^{2}Q_{T}\left(
\beta _{T}^{0}\right)}{\partial \beta \partial \beta ^{\prime }} \sqrt{T}\left( \hat{\beta}_{T}^{r}-\beta
_{T}^{0}\right) +\frac{\partial g^{\prime }(\beta _{T}^{0})}{\partial \beta }%
\sqrt{T}\hat{\lambda}_{T} &=&o_{P}(1).  \label{FOC} 
\end{eqnarray}
Recalling the definition of the infeasible unconstrained estimator $\ddot{%
	\beta}_{T}$, we can rewrite the LHS of the above equation as follows (while
the RHS is exactly zero when the function $Q_{T}$ is quadratic and the
constraints are linear): 
\[
J_{T}\sqrt{T}\left( \ddot{\beta}_{T}-\beta _{T}^{0}\right) -J_{T}\sqrt{T}%
\left( \hat{\beta}_{T}^{r}-\beta _{T}^{0}\right) +\frac{\partial {g}%
_{T}^{\prime }(\beta _{T}^{0})}{\partial \beta }\sqrt{T}\hat{\lambda}%
_{T}=o_{P}(1) .
\]
By Lemma 1, and Assumption \textbf{A1}, all
three terms of the LHS of the above equality are all $O_{P}(1).$ We deduce that
\begin{equation}
J_{T}\sqrt{T}\left( \hat{\beta}_{T}^{r}-\beta _{T}^{0}\right) -\frac{%
\partial g^{\prime }(\beta _{T}^{0})}{\partial \beta }\sqrt{T}\hat{\lambda}%
_{T}=J_{T}\sqrt{T}\left( \ddot{\beta}_{T}-\beta _{T}^{0}\right) +o_{P}(1).
\label{decomp}
\end{equation}Moreover, as already noted above, the remainder term $o_{P}(1)$\ in (\ref{decomp})
is zero when the criterion function $Q_{T}$ is quadratic and the constraints are linear. 
\end{proof}

\begin{proof}[Proof of Theorem 1]
By definition 
\[
\sqrt{T}\frac{\partial Q_{T}(\hat{\beta}_{T}^{r})}{\partial \beta ^{ }}%
+\frac{\partial ^{2}Q_{T}(\hat{\beta}_{T}^{r})}{\partial \beta \partial
\beta ^{\prime }}\sqrt{T}\left( \widehat{\beta }_{T}-\hat{\beta}%
_{T}^{r}\right) =0
\]%
Therefore, by a Taylor expansion of the first term around the true value $%
\beta _{T}^{0}:$%
\[
\sqrt{T}\frac{\partial Q_{T}(\beta _{T}^{0})}{\partial \beta ^{ }}+%
\frac{\partial ^{2}Q_{T}(\beta _{T}^{0})}{\partial \beta \partial \beta
^{\prime }}\sqrt{T}(\hat{\beta}_{T}^{r}-\beta _{T}^{0})+\frac{\partial
^{2}Q_{T}(\hat{\beta}_{T}^{r})}{\partial \beta \partial \beta ^{\prime }}%
\sqrt{T}\left( \widehat{\beta }_{T}-\hat{\beta}^{r}_T\right) = o_{P}(1)
\]%
and then, since $\sqrt{T}(\hat{\beta}_{T}^{r}-\beta _{T}^{0})=O_{P}(1)$, we
can obviously simplify the above decomposition to obtain
\[
\sqrt{T}\frac{\partial Q_{T}(\beta _{T}^{0})}{\partial \beta ^{ }}=-%
\frac{\partial ^{2}Q_{T}(\hat{\beta}_{T}^{r})}{\partial \beta \partial \beta
^{\prime }}\sqrt{T}\left( \widehat{\beta }_{T}-\beta _{T}^{0}\right)
+o_{P}(1).
\]%
Since by Assumption \textbf{A1}, we know that%
\[
\plim_{T\rightarrow\infty }\frac{\partial ^{2}Q_{T}(\hat{\beta}_{T}^{r})}{\partial
\beta \partial \beta ^{\prime }}=-\mathcal{J}^{0}
\]%
we can conclude that $\sqrt{T}\left( \widehat{\beta }_{T}-\beta
_{T}^{0}\right) =O_{P}(1)$ and%
\[
\sqrt{T}\left( \widehat{\beta }_{T}-\beta _{T}^{0}\right) =\left[ \mathcal{J}%
^{0}\right] ^{-1}\sqrt{T}\frac{\partial Q_{T}(\beta _{T}^{0})}{\partial
\beta ^{ }}+o_{P}(1)
\]%
By comparison with the definition of $\ddot{\beta}_{T}$:%
\[
\sqrt{T}\left( \ddot{\beta}_{T}-\beta _{T}^{0}\right) =\left[ J_{T}\right]
^{-1}\sqrt{T}\frac{\partial Q_{T}(\beta _{T}^{0})}{\partial \beta ^{ }}
\]%
we have the announced equivalence between estimators. 
\end{proof}

\begin{proof}[Proof of Proposition 2]
For $\beta _{T}^{\ast }$ a component-by-component intermediate value between $\beta^0_T$ and $\hat{\beta}^r_T$, by Assumption \textbf{A1}(ii), we
deduce that
\[
\left\Vert \frac{\partial ^{2}Q_{T}(\beta _{T}^{\ast })}{\partial \beta
	\partial \beta ^{\prime }}\left[ \frac{\partial ^{2}Q_{T}(\hat{\beta}%
	_{T}^{r})}{\partial \beta \partial \beta ^{\prime }}\right]
^{-1}-\text{Id}_{d_{\beta }}\right\Vert =o_{P}(1),
\]where $\text{Id}_{d_{\beta }}$ denotes the $d_\beta\times d_\beta$ identity matrix. Moreover, this bound does not depend on $\theta $. We apply the fact
that this bound remains uniformly valid on $\Theta $\ when quantities are
multiplied by continuous functions of $\theta $\ over compact $\Theta $.

For all $\theta \in \Theta $,
\begin{eqnarray*}
	\sqrt{T}\bar{m}_{TH}[\theta ;\widehat{\beta }_{T}] &=&\sqrt{T}\frac{\partial
		Q_{TH}(\theta ,\beta _{T}^{0})}{\partial \beta }-\frac{\partial
		^{2}Q_{TH}(\theta ,\hat{\beta}_{T}^{r})}{\partial \beta \partial \beta
		^{\prime }}\left[ \frac{\partial ^{2}Q_{T}(\hat{\beta}_{T}^{r})}{\partial
		\beta \partial \beta ^{\prime }}\right] ^{-1}\sqrt{T}\frac{\partial
		Q_{T}(\beta _{T}^{0})}{\partial \beta } \\
	&&+\left\{ \frac{\partial ^{2}Q_{TH}(\theta ,\beta _{T}(\theta ))}{\partial
		\beta \partial \beta ^{\prime }}-\frac{\partial ^{2}Q_{TH}(\theta ,\hat{\beta%
		}_{T}^{r})}{\partial \beta \partial \beta ^{\prime }}\right\} \sqrt{T}\left(	\hat{\beta}_{T}^{r}-\beta _{T}^{0}\right) +o_{P}(1),
\end{eqnarray*}
where $\beta _{T}(\theta )$\ is a component-by-component intermediate
value between $\beta _{T}^{0}$\ and $\hat{\beta}_{T}^{r}$, and where the bound 
$o_{P}(1)$\ does not depend on $\theta $. By Assumption \textbf{A3}(iii)
\[
\sup_{\theta \in \Theta }\left\Vert \frac{\partial ^{2}Q_{TH}(\theta ,\beta
	_{T}(\theta ))}{\partial \beta \partial \beta ^{\prime }}-\frac{\partial
	^{2}Q_{TH}(\theta ,\hat{\beta}_{T}^{r})}{\partial \beta \partial \beta
	^{\prime }}\right\Vert =o_{P}(1), 
\]so that, by Assumption \textbf{A3}(iii) and Assumption \textbf{A1}(ii), we have
\[
\sup_{\theta \in \Theta }\left\Vert \frac{\partial ^{2}Q_{TH}(\theta ,\hat{%
		\beta}_{T}^{r})}{\partial \beta \partial \beta ^{\prime }}\left[ \frac{%
	\partial ^{2}Q_{T}(\hat{\beta}_{T}^{r})}{\partial \beta \partial \beta
	^{\prime }}\right] ^{-1}-\text{Id}_{d_{\beta }}\right\Vert =o_{P}(1) .
\]
Therefore, we can conclude that
\[
\sup_{\theta \in \Theta }\left\Vert \sqrt{T}\bar{m}_{TH}[\theta ;\widehat{%
	\beta }_{T}]-\left\{ \sqrt{T}\frac{\partial Q_{TH}(\theta ,\beta _{T}^{0})}{%
	\partial \beta }-\sqrt{T}\frac{\partial Q_{T}(\beta _{T}^{0})}{\partial
	\beta }\right\} \right\Vert =o_{P}(1) .
\]A similar argument would allow us to prove
\[
\sup_{\theta \in \Theta }\left\Vert \sqrt{T}m_{TH}^{CFS}[\theta ;\hat{\lambda%
}_{T}]-\left\{ \sqrt{T}\frac{\partial Q_{TH}(\theta ,\beta _{T}^{0})}{%
	\partial \beta }-\sqrt{T}\frac{\partial Q_{T}(\beta _{T}^{0})}{\partial
	\beta }\right\} \right\Vert =o_{P}(1). 
\]

Now, revisiting the definition of $\partial Q_{TH}(\theta,\ddot{\beta}_{T})/\partial\beta$, for all $\theta \in \Theta$,
\begin{eqnarray*}
	\sqrt{T}\frac{\partial Q_{TH}(\theta ,\ddot{\beta}_{T})}{\partial \beta } &=&%
	\sqrt{T}\frac{\partial Q_{TH}(\theta ,\beta _{T}^{0})}{\partial \beta }+%
	\frac{\partial ^{2}Q_{TH}(\theta ,\beta _{T}(\theta ))}{\partial \beta
		\partial \beta ^{\prime }}\sqrt{T}\left( \ddot{\beta}_{T}-\beta _{T}^{0}%
	\right) \\
	&=&\sqrt{T}\frac{\partial Q_{TH}(\theta ,\beta _{T}^{0})}{\partial \beta }+%
	\frac{\partial ^{2}Q_{TH}(\theta ,\beta _{T}(\theta ))}{\partial \beta
		\partial \beta ^{\prime }}\left\{ -\frac{\partial ^{2}Q_{T}(\beta _{T}^{0})}{%
		\partial \beta \partial \beta ^{\prime }}\right\} ^{-1}\sqrt{T}\frac{%
		\partial Q_{T}(\beta _{T}^{0})}{\partial \beta },
\end{eqnarray*}where $\beta _{T}(\theta )$\ is a component-by-component intermediate
value between $\beta _{T}^{0}$\ and $\ddot{\beta}_{T}$.
Applying a similar argument to the one above, and using the fact that $%
\sqrt{T}( \ddot{\beta}_{T}-\beta _{T}^{0}) =O_{P}(1)$, we have 
\[
\sup_{\theta \in \Theta }\left\Vert \sqrt{T}\frac{\partial Q_{TH}(\theta ,%
	\ddot{\beta}_{T})}{\partial \beta }-\left\{ \sqrt{T}\frac{\partial
	Q_{TH}(\theta ,\beta _{T}^{0})}{\partial \beta }-\sqrt{T}\frac{\partial
	Q_{T}(\beta _{T}^{0})}{\partial \beta }\right\} \right\Vert =o_{P}(1) .
\]

Therefore, the three estimating equations $\sqrt{T}\bar{m}_{TH}[\theta ;\widehat{\beta }_{T}],%
\sqrt{T}m_{TH}^{CFS}[\theta ;\hat{\lambda}_{T}]$ and $\sqrt{T}\frac{\partial
	Q_{TH}(\theta ,\ddot{\beta}_{T})}{\partial \beta }$\ are each asymptotically
equivalent to $\left\{ \sqrt{T}\frac{\partial
	Q_{TH}(\theta ,\beta _{T}^{0})}{\partial \beta }-\sqrt{T}\frac{\partial
	Q_{T}(\beta _{T}^{0})}{\partial \beta }\right\} $, and the result follows.
\end{proof}

\begin{proof}[Proof of Proposition 3]
(i) We first prove that $\widehat{\theta }_{T,H}^{CFS}(W)$\ is consistent.
By Assumption \textbf{A4}(ii),  $m_{TH}^{CFS}[\theta ;\hat{\lambda}_{T}]$ converges in
probability, uniformly on $\theta \in \Theta $, towards
\begin{flalign*}
\plim_{T\rightarrow\infty }\left\{\frac{\partial Q_{TH}(\theta ,\hat{\beta}_{T}^{r})}{%
\partial \beta }-\frac{\partial Q_{T}(\hat{\beta}_{T}^{r})}{%
\partial \beta }\right\}&=L(\theta ,\beta ^{0})-L(\theta ^{0},\beta ^{0})=L(\theta ,\beta ^{0})
\end{flalign*}
The identification Assumption \textbf{A4}(iii), jointly with compactness of $\Theta $
and the continuity assumption \textbf{A4}(i), then yields
\[
\Plim{T\rightarrow\infty }\left\{\widehat{\theta }_{T,H}^{CFS}(W)\right\}=\theta ^{0}.
\]

\noindent(ii) By comparing \eqref{scoreus} and \eqref{CFS},  we have 
\begin{flalign*}
\bar{m}_{TH}[\theta ;\widehat{\beta }_{T}]-m_{TH}^{CFS}[\theta ;\hat{%
\lambda}_{T}]&=\frac{\partial ^{2}Q_{TH}(\theta ,\hat{\beta}_{T}^{r})}{\partial \beta
\partial \beta ^{\prime }}\left( \widehat{\beta }_{T}-\hat{\beta}%
_{T}^{r}\right) +\frac{\partial Q_{T}(\hat{\beta}_{T}^{r})}{\partial \beta }
\\
&=\left\{\text{Id}_{d_{\beta}}-\frac{\partial ^{2}Q_{TH}(\theta ,\hat{\beta}_{T}^{r})}{\partial \beta
\partial \beta ^{\prime }}\left[ \frac{\partial ^{2}Q_{T}(\hat{\beta}%
_{T}^{r})}{\partial \beta \partial \beta ^{\prime }}\right] ^{-1}\right\}\frac{%
\partial Q_{T}(\hat{\beta}_{T}^{r})}{\partial \beta },
\end{flalign*}By {Assumptions} \textbf{A1}(ii) and \textbf{A3}(iii), this difference converges, uniformly on $\theta
\in \Theta $, towards
\[
-\mathcal{J}(\theta ,\beta ^{0})[\mathcal{J}^{0}]^{-1}L(\theta ^{0},\beta
^{0})+L(\theta ^{0},\beta ^{0})=0
\]
{where} $\mathcal{J}^{0}=\mathcal{J}(\theta^0,\beta^0)$. Then, by a standard argument (see, e.g., Pakes and Pollard, 1989,
page 1038), we deduce that
\[
\Plim{T\rightarrow\infty }\left\{\widehat{\theta }_{T,H}^{s}(W)\right\}=\Plim{T\rightarrow\infty }\left\{\widehat{%
\theta }_{T,H}^{CFS}(W)\right\}=\theta ^{0}.
\]
\noindent(iii) By Assumptions \textbf{A1}(ii) and \textbf{A3}(iii), 
\[
\sup_{\left\Vert \theta -\theta ^{0}\right\Vert \leq \gamma /\sqrt{T}}\left\Vert -%
\frac{\partial ^{2}Q_{TH}(\theta ,\hat{\beta}_{T}^{r})}{\partial \beta
\partial \beta ^{\prime }}\left[ \frac{\partial ^{2}Q_{T}(\hat{\beta}%
_{T}^{r})}{\partial \beta \partial \beta ^{\prime }}\right]
^{-1}+\text{Id}_{d_{\beta }}\right\Vert =o_{P}(1).
\]
Then, deduce from the above decomposition that
\[
\sup_{\left\Vert \theta -\theta ^{0}\right\Vert \leq \gamma/\sqrt{T}}\left\Vert 
\bar{m}_{TH}[\theta ;\widehat{\beta }_{T}]-m_{TH}^{CFS}[\theta ;\hat{\lambda}%
_{T}]\right\Vert =o_{P}\left( \left\Vert \frac{\partial Q_{T}(\hat{\beta}%
_{T}^{r})}{\partial \beta }\right\Vert \right) =o_{P}\left( \frac{1}{\sqrt{T}%
}\right) ,
\]
which in turn implies that
\[
\sup_{\left\Vert \theta -\theta ^{0}\right\Vert \leq \gamma/\sqrt{T}}\left\vert
S_{T}^{unr}(\theta )-S_{T}^{res}(\theta )\right\vert =o_{P}(1/{T}),
\]
for $S_{T}^{unr}(\theta )$\ and $S_{T}^{res}(\theta )$ respectively
the objective functions minimized in (\ref{unr}) and (\ref{res}) to define
the estimators $\widehat{\theta }_{T,H}^{s}(W)$\ and $\widehat{\theta }%
_{T,H}^{CFS}(W)$\ respectively.

It is then a standard argument (see, e.g., Pakes and Pollard, 1989, page 1040)
to deduce that, using the asymptotic normality in Assumption \textbf{A1}(iii), the corresponding extremum estimators are asymptotically
equivalent: $
\left\Vert \widehat{\theta }_{T,H}^{s}(W)-\widehat{\theta }%
_{T,H}^{CFS}(W)\right\Vert =o_{P}(1/\sqrt{T}).$
\end{proof}

\begin{proof}[Proof of Theorem 2]
We first show that $\sqrt{T}\left(\widehat{\theta}^s_{T,H}(W)-\theta_T\right)=O_{P}(1)$. Let $\widehat{\theta}_T:=\widehat{\theta}^s_{T,H}(W)$ and define, for a vector $x$, $\|x\|_{W_{}}=\sqrt{x'Wx}$. By the triangle inequality, 
\begin{flalign}\label{eqsnew1}
\|L(\widehat{\theta}_{T},\beta^0_T)-L(\theta_T,\beta^0_T)\|_{W_{}}&\leq \left\|L(\widehat{\theta}_{T},\beta^0_T)-L(\theta_T,\beta^0_T)-\bar{m}_{TH}[\widehat{\theta}_{T};\widehat{\beta}_{T}] \right\|_{W_{}}+\|\bar{m}_{TH}[\widehat{\theta}_{T};\widehat{\beta}_{T}]\|_{W_{}}.
\end{flalign}Consider the first term in \eqref{eqsnew1}. Recall the definition of $\bar{m}_{TH}[\widehat{\theta}_{T};\widehat{\beta}_{T}]$, and apply the definitions of $\widehat{\beta}_{T},\hat{\beta}^r_{T}$, and the uniform convergence in Assumption \textbf{A3}(iii), to deduce 
\begin{flalign}
\bar{m}_{TH}[\widehat{\theta}_{T};\widehat{\beta}_{T}]&=\frac{\partial Q_{TH}(\widehat{\theta}_{T},\hat{\beta}^r_{T})}{\partial\beta}+\frac { \partial ^ { 2 } Q _ { T H } \left( \widehat{\theta}_{T} , \hat{\beta}^{r}_{T}\right) } { \partial \beta \partial \beta ^ { \prime } }\left(\widehat{\beta}_{T}-\hat{\beta}^r_{T}\right)
\nonumber\\&=\frac{\partial Q_{TH}(\widehat{\theta}_{T},\hat{\beta}^r_{T})}{\partial\beta}-\mathcal{J}(\theta_T,\beta^0_T)[\mathcal{J}^0]^{-1}\frac{\partial Q_{T}(\hat{\beta}^r_T)}{\partial\beta}+o_{P}(1/\sqrt{T}).\label{eqsnew2}
\end{flalign}Applying equation \eqref{eqsnew2} and the triangle inequality, we obtain
\begin{flalign}
\left\|L(\widehat{\theta}_{T},\beta^0_T)-L(\theta_T,\beta^0_T)-\bar{m}_{TH}[\widehat{\theta}_{T};\widehat{\beta}_{T}] \right\|_{W_{}}&\leq\sup_{\theta\in\Theta,\|\beta-\beta^0_{T}\|\leq\frac{\gamma}{\sqrt{T}}}\left\|\frac{\partial Q_{TH}({\theta},{\beta})}{\partial\beta} -L({\theta},\beta^0_T)\right\|_{W_{}}\nonumber\\&+\left\|\frac{\partial Q_{T}(\hat{\beta}^r_T)}{\partial\beta}-L(\theta_T,\beta^0_T)\right\|_{W}\label{eqsnew3}.
\end{flalign}In \eqref{eqsnew3}, the first term is $O_{P}(1/\sqrt{T})$ by \textbf{A4}(ii), and the second term is $O_{P}(1/\sqrt{T})$ by Assumption \textbf{A1}(iii) and \textbf{Lemma 1}. Analyzing the second term in equation \eqref{eqsnew1}, note that, by definition $$\|\bar{m}_{TH}[\widehat{\theta}_{T};\widehat{\beta}_{T}]\|_{W_{}}\leq \|\bar{m}_{TH}[\theta_T;\widehat{\beta}_{T}]\|_{W_{}}.$$Apply the same decomposition in \eqref{eqsnew2} to the term $\bar{m}_{TH}[\theta_T;\widehat{\beta}_{T}]$, and the triangle inequality to obtain 
$$\|\bar{m}_{TH}[\widehat{\theta}_{T};\widehat{\beta}_{T}]\|_{W_{}}\leq \|\bar{m}_{TH}[\theta_T;\widehat{\beta}_{T}]\|_{W_{}}\leq \left\|\frac{\partial Q_{TH}({\theta}^0_{T},{\beta}^0_{T})}{\partial\beta}\right\|_{W}+O_{P}(1/\sqrt{T}).$$ From Assumption \textbf{A3}(ii), we  have $\|{\partial Q_{TH}({\theta}^0_{T},{\beta}^0_{T})}/{\partial\beta}\|_{}=O_{P}(1/\sqrt{T})$, which yields 
 \begin{equation}\label{eqsnew4}
 \|\bar{m}_{TH}[\widehat{\theta}_{T};\widehat{\beta}_{T}]\|_{W_{}}\leq O_{P}(1/\sqrt{T}).
 \end{equation}Applying the results in \eqref{eqsnew3} and \eqref{eqsnew4} into \eqref{eqsnew2}, we arrive at $\|L(\widehat{\theta}_{T},\beta^0_T)-L({\theta}^{},\beta^0_T)\|_{W_{}}=O_{P}(1/\sqrt{T})$. From the local identification Assumption \textbf{A5}, we then have that, for some constant $C>0$, $$C\|\widehat{\theta}_{T}-\theta_T\|\leq\|L(\widehat{\theta}_{T},\beta^0_T)-L(\theta_T,\beta^0_T)\|_{W_{}}=O_{P}(1/\sqrt{T}).$$

Having proven $\sqrt{T}\left(\widehat{\theta}_{T}-\theta_T\right)=O_{P}(1/\sqrt{T})$, the remainder of the proof proceeds through a standard first-order Taylor series of the first-order conditions. First, apply the result of \textbf{Proposition 2} to obtain
\begin{flalign*}
0&=\left[\frac{\partial^2 Q_{TH}(\widehat{\theta}_{T},\hat{\beta}^r_T)}{\partial\beta\partial\theta'}\right]'W\sqrt{T}\bar{m}_{TH}[\widehat{\theta}_{T};\widehat{\beta}_{T}]+o_{P}(1)=\left[\frac{\partial^2 Q_{TH}(\widehat{\theta}_{T},\hat{\beta}^r_T)}{\partial\beta\partial\theta'}\right]'W\sqrt{T}\frac{\partial Q_{TH}(\widehat{\theta}_{T},\ddot{\beta}_T)}{\partial\beta}+o_{P}(1).
\end{flalign*}
Now, disregarding terms of smaller order than $O_{P}(1/\sqrt{T})$, a Taylor series expansion of $\sqrt{T}{\partial Q_{TH}(\widehat{\theta}_{T},\ddot{\beta}_T)}/{\partial\beta}$ around $\theta_T$ yields
\begin{flalign*}
0&= \left[\frac{\partial L(\theta_T,\beta^0_T) }{\partial\theta'}\right]'W\sqrt{T}\left\{\frac{\partial Q_{TH}(\theta_T,\ddot{\beta}_{T})}{\partial\beta}+\frac{\partial Q_{TH}(\theta_T,\ddot{\beta}_{T})}{\partial\beta\partial\theta'}\left(\widehat{\theta}_{T}-\theta_T\right)\right\}+o_{P}(1).
\end{flalign*}Expanding $\sqrt{T}\partial Q_{TH}(\theta_T,\ddot{\beta}_T)/{\partial\beta}$ as in equation \eqref{final}, and making use of the uniform convergence in Assumption \textbf{A3}(iii), we obtain (up to an $o_{P}(1)$ term)
\begin{flalign*}
0&=\left[\frac{\partial L(\theta_T,\beta^0_T) }{\partial\theta'}\right]'W\sqrt{T}\left\{\frac{\partial Q_{TH}(\theta_T,{\beta}^0_{T})}{\partial\beta}-\frac{\partial Q_{T}({\beta}^0_{T})}{\partial\beta}\right\}+\left[\frac{\partial L(\theta_T,\beta^0_T) }{\partial\theta'}\right]'W\left[\frac{\partial L(\theta_T,\beta^0_T) }{\partial\theta'}\right]\sqrt{T}\left(\widehat{\theta}_{T}-\theta_T\right)
\end{flalign*}Rearranging terms and making use of Assumption \textbf{A5}, 
\begin{flalign*}
\sqrt{T}\left(\widehat{\theta}_{T}-\theta_T\right)=-\left\{\left[\frac{\partial L(\theta_T,\beta^0_T) }{\partial\theta'}\right]'W\left[\frac{\partial L(\theta_T,\beta^0_T) }{\partial\theta'}\right]\right\}^{-1}\left[\frac{\partial L(\theta_T,\beta^0_T) }{\partial\theta'}\right]'W\sqrt{T}\left\{\frac{\partial Q_{TH}(\theta_T,{\beta}^0_{T})}{\partial\beta}-\frac{\partial Q_{T}({\beta}^0_{T})}{\partial\beta}\right\}.
\end{flalign*}By Assumptions \textbf{A1}(iii) and \textbf{A3}(ii), the term in brackets is an asymptotically Gaussian mean-zero random variable. The stated result then follows. 
\end{proof}

\begin{proof}[Proof of Theorem 3]
The result follows from the following sequence of arguments:  {(i)} $\widehat{\theta}^{c}_{T,H}$ solves 
$\widehat{\beta}_{T}=\tilde{\beta}^{c}_{TH}(\theta)$; \textbf{(ii)} $\widehat{\theta}^{s}_{T,H}$ solves 
$0=\bar{m}_{TH}[\theta,\widehat{\beta}_{T}]$; \textbf{(iii)} From \textbf{(ii)} and the structure of $\bar{m}_{TH}[\theta,\widehat{\beta}_{T}]$ we have, re-arranging $0=\bar{m}_{TH}[\widehat{\theta}^{s}_{T},
\widehat{\beta}_{T}]$ and solving for $\widehat{\beta}_{T}$, $$\widehat{\beta}_{T}=\hat{\beta}^{r}_{T}-\left[\frac{\partial^{2}Q_{TH}[\widehat{\theta}^{s}_{T},\hat{\beta}^{r}_{T}]}{\partial\beta\partial\beta^{\prime}}\right]^{-1}\frac{\partial^{}Q_{TH}[\widehat{\theta}^{s}_{T},\hat{\beta}^{r}_{T}]}{\partial\beta}= \tilde{\beta}^{c}_{TH}(\widehat{\theta}_{T}^{s}),$$ 
where the last equality follows from the definition of $\tilde{\beta}^{c}_{TH}(\theta)$.  Therefore, from \textbf{(i)} we have $ \widehat{\beta}_{T}=\tilde{\beta}^{c}_{TH}(\widehat{\theta}^{c}_{T})$ and from \textbf{(iii)} we have 
$
\widehat{\beta}_{T}=\tilde{\beta}^{c}_{TH}(\widehat{\theta}^{c}_{T})=\tilde{\beta}^{c}_{TH}(\widehat{\theta}_{T}^{s}).
$\end{proof}

\end{document}